%
\input amstex
\documentstyle {amsppt}
\loadbold
\font\smbf=cmbx10 at 11pt
\font\medbf=cmbx10 scaled \magstep1
\font\bigbf=cmbx10 scaled \magstep2
\magnification=1000
\TagsOnRight
\input psfig.sty
\hsize=6in
\vsize=8in
\hfuzz=3 pt
\hoffset=0.1in
\tolerance 9000
\def\L{\text{{\bf L}}}
\def\meas{\hbox{meas}}
\def\vvs{\vskip 0.5em}
\def\v{\vskip 1em}
\def\vs{\vskip 2em}
\def\vsk{\vskip 4em}
\def\a{\alpha}
\def\tv{\hbox{\text{Tot.Var.}}}

\def\bn{\text{\bf n}}
\def\la{\big\langle}

\def\ra{\big\rangle}
\def\ov{\overline}

\def\R{\Bbb R}
\def\A{{\Cal A}}
\def\O{{\Cal O}}

\def\BV{{\strut \text{{\rm BV}}}}
\def\ds{\diamondsuit}
\def\di{\sharp}
\def\bem{\flat}
\def\beq{\natural}

\def\Om{\Omega}

\def\forallt{\hbox{for all~}}

\def\ve{\varepsilon}

\def\n{\noindent}

\def\fine{\hfill $\square$}

%
\voffset=0.4in
%
%
\def\ABpv{\hbox{[1]}}
%
\def\ABfs{\hbox{[2]}}
%

%

%
%
%
\def\Bimp{\hbox{[3]}}
%
%

%
\def\CLSS{\hbox{[5]}}
%

%
%
\def\CLSWnac{\hbox{[6]}}
%
%
\def\CLRSfslf{\hbox{[4]}}
%

%
%

%

%

%

%

%

%
%

%
%

%
%

%
%

%
%

%
%

%
%

%
%
%
\def\Riscl{\hbox{[7]}}
%
%
%
\def\Risclsf{\hbox{[8]}}
%

%
%

%

%

%
\def\Ss{\hbox{[9]}}
%

%


\def\rhead=\vbox to 2truecm{\hbox
{\ifodd\pageno\rightheadline\else\leftheadline\fi}
\hbox to 1in{\hfil}\psill}
\tenrm
\def\rightheadline{\hfil{\eightrm Stability Rates for Patchy Vector Fields
}\hfil\folio}
\def\leftheadline{\folio\hfil{\eightrm F. Ancona and A. Bressan}\hfil}
\def\onepageout#1{\shipout\vbox
\offinterlineskip
\vbox to 2in{\rhead}
\vbox to \pageheight
\advancepageno}
\NoBlackBoxes
%
\document
\topskip=20pt
\nopagenumbers

\null
\v
\vskip 2cm
\centerline{\bigbf Stability Rates for Patchy Vector Fields}
\vskip 2cm 
\centerline{\it Fabio Ancona$^{(*)}$ and Alberto Bressan$^{(**)}$}
\vs

\parindent 40pt

\item{(*)} Dipartimento di Matematica and C.I.R.A.M., Universit\`a di Bologna,
\item{} Piazza Porta S.~Donato~5,~Bologna~40127,~Italy.
\item{} e-mail: ancona\@ciram3.ing.unibo.it.
\v
\item{(**)} S.I.S.S.A., Via Beirut~4,~Trieste~34014,~Italy.
\item{} e-mail: bressan\@sissa.it.
\vskip 1cm

\centerline{June 2001}
\vskip 2cm

\n {\bf Abstract.} 
The paper is concerned with the stability of the set of trajectories
of a {\it patchy} vector field, in the presence of impulsive
perturbations.  
Patchy vector fields are
discontinuous, piecewise smooth vector fields
that were introduced in \ABpv\ 
to study feedback stabilization problems.
For patchy vector fields in the plane, with polygonal
patches in generic position, we show that the distance 
between a perturbed trajectory and an unperturbed one 
is of the same order of magnitude as the impulsive forcing term.
\vs
\vs

\noindent
{\it 1991 AMS-Subject Classification.} \
%
%
34 A,
%
%
34 D, 
%
%
49 E, 
%
%
93 D.
\vfill\eject
\vsk

\pageno=1

\parindent 20pt

\n{\medbf 1 - Introduction.}
\v
Let $g$ be a bounded vector field, and consider the
Cauchy problem with impulsive perturbations
$$
\dot y = g(y)+ \dot w.
\tag 1.1
$$
Here $w=w(t)$ is a
left continuous function with bounded variation.
By a solution of (1.1) 
with initial condition
$$
y(t_0)=y_0,
\tag 1.2
$$
we mean
a measurable function $t\mapsto y(t)$ such that
$$
y(t)=y_0+\int_{t_0}^t g\big(y(s)\big)\,ds +
\big[w(t)-w(t_0)\big].
\tag 1.3
$$
If $w(\cdot)$ is discontinuous, the forcing term in (1.1) will have
impulsive behavior, and the solution $y(\cdot)$
will be discontinuous as well.
We choose to work with (1.1) because it provides a simple and general
framework to study stability properties.  
Indeed, consider a system with
both inner and outer perturbations, of the form
$$
\dot x=g\big( x + e_1(t)\big) + e_2(t).
\tag 1.4
$$
Then, the map $y=y(t)\doteq x(t)+e_1(t)$ satisfies 
the impulsive equation
$$
\dot y=g(y)+e_2(t)+\dot e_1(t)=g(y)+\dot  w,
$$
where
$$ 
w (t)=e_1(t)+\int_{t_0}^t e_2(s)\,ds.
$$
Therefore, from the 
stability of solutions of (1.1) under small BV
perturbations $w$, one can immediately \, deduce a 
result on the stability of solutions of (1.4), when 
\, $\tv\{e_1\}$ \,
and \, $\|e_2\|_{\L^1}$ \, are suitably 
small.

Our main concern is how much a trajectory is
affected by the presence of the impulsive perturbation.
More precisely, we wish to estimate the distance, in the $\L^\infty$
norm, between solutions of the two Cauchy problems
$$
\left\{\eqalign{\dot x&=g(x),\cr x(0)&=x_0,\cr}\right.\qquad\qquad
\left\{\eqalign{\dot y&=g(y)+\dot w,\cr y(0)&=x_0,\cr}\right.
\eqno(1.5)
$$
Consider first the special case where $g$ is a 
continuous vector field with
Lipschitz constant $L$.  It is then well known that
the Cauchy problems (1.5)
have unique solutions, obtained by a fixed point argument (see \Bimp).
Their distance can be estimated as
$$
\big| y(t)-x(t)\big|\leq \int_0^t  e^{L(t-s)}\,\big| dw(s)\big|
\leq e^{Lt}\cdot \tv\{w\}.
\tag 1.6
$$
In other words, on a fixed time interval, this distance grows
linearly with $\tv\{w\}$.

In this paper, we will prove a similar estimate 
in the case where
$g$ is a discontinuous, {\it patchy vector field}.
These vector fields were introduced in \ABpv\ in order to study 
feedback stabilization problems.  We recall the main definitions:
\v
\n{\bf Definition 1.1.} By a {\it patch} we mean a pair 
$\big(\Omega,\, g\big)$ where 
$\Omega\subset\R^n$ 
is an open domain with smooth boundary $\partial\Om,$ and 
$g$ is a smooth vector field 
defined on a neighborhood of the closure
$\ov\Omega$, which points strictly inward at each
boundary point $x\in\partial\Om$.
\v  
Calling $\bn(x)$ the outer normal at the boundary
point $x$, we thus require 
$$
\la g(x),~\bn(x)\ra <0\qquad\forallt\quad x\in\partial\Om.
\tag 1.7
$$
\v
\n{\bf Definition 1.2.}  \ We say that $g:\Omega\mapsto\R^n$
is a {\it patchy vector field} on the open domain $\Om$
if there exists a family of patches
$\big\{ (\Omega_\alpha,~g_\alpha) ;~~ \alpha\in\A\big\}$ such that

\n - $\A$ is a totally ordered set of indices, 

\n - the open sets $\Omega_\alpha$ form a locally finite covering of 
$\Omega$, i.e. $\Omega=\cup_{\alpha\in \A}\Omega_\alpha$
and every compact set~$K\subset \Bbb R^n$
intersect only a finite number of domains $\Omega_\alpha, \, \alpha \in \A$;

\n - the vector field $g$ can be written in the form
$$
g(x) = g_\alpha (x)\qquad \hbox{if}\qquad x \in 
\Omega_\alpha \setminus 
\displaystyle{\bigcup_{\beta > \alpha} \Omega_\beta}.
\tag 1.8              
$$
\v
By setting
$$
\alpha^*(x) \doteq \max\big\{\alpha \in \A~;~~ x \in \Omega_\alpha
\big\},
\tag 1.9
$$
we can write (1.8) in the equivalent form
$$
g(x) = g_{_{\alpha^*(x)}}(x) \qquad \forallt~~x \in \Omega.
\tag 1.10
$$
We shall occasionally adopt the longer notation
$\big(\Omega,\ g,\ (\Omega_\alpha,\,g_\alpha)_{_{\alpha\in \A}} \big)$ 
to indicate a patchy vector field, specifying both the domain
and the single patches.
If $g$ is a patchy vector field, the differential equation
$$
\dot x = g(x)
\tag 1.11
$$
has many interesting properties.
In particular, in \ABpv\ it was proved that the set 
of Carath\'eodory solutions of
(1.11) is closed (in the topology of uniform convergence) 
but possibly not
connected.  
Moreover, given an initial condition
$
x(t_0)=x_0,
$
the corresponding Cauchy problem has at least one forward
solution, and at most one backward solution, in the Carath\'eodory
sense. For every  Carath\'eodory solution \, $x(\cdot)$ \, of (1.11),
the map \, $t \mapsto \a^*(x(t))$ \, is left continuous
and non-decreasing.

Since the Cauchy problem for (1.11) does not have 
forward uniqueness and continuous dependence, one clearly cannot
expect that a single solution can be stable under small
perturbations.
Instead, one can establish 
the following stability  property referring to the whole set
of solutions.
\vs
\n{\bf Proposition 1.3} [2, \, Corollary~1.3] \
{\it Let $g$ be a 
patchy vector field
on an open domain $\Omega\subset \R^n.$
Given any closed subset $A\subset \Omega,$ 
any compact set $K\subset A,$ and every $T, \ve>0$, 
there exists \,
$\delta>0$ \, such that the following holds.
If $y:[0,T]\mapsto A$ is a solution of the perturbed system (1.1)
with $y(0)\in K$ and $\tv (w)<\delta$, then there exists a solution
$x:[0,T]\mapsto \Omega$ of the unperturbed equation (1.11) with}
$$
\big\|x-y\big\|_{\L^\infty([0,T])}
<\ve\,.\tag 1.12
$$
The relevance of this result for the robustness of discontinuous
feedback controls is discussed in \ABfs.
\v
In connection with Proposition 1.3, it is interesting to study how
the distance $\|y-x\|_{\L^\infty}$ can depend on the
perturbation $w$.  For a general BV function $w$,
the derivative $\dot w$ is a Radon measure whose total mass coincides
with the total variation of $w$. It is thus natural to use the BV norm
$\|w\|_{BV}$ as a measure of the strength of the perturbation.
In the case of a Lipschitz continuous field $g$, 
we have seen in 
(1.6) that this distance grows linearly with $\|w\|_{BV}$\,.
In the case of patchy vector fields, one cannot expect a linear
dependence, in general.
\v
\n{\bf Example 1.4.}  Consider a patchy vector field on $\R^2$,
as in figure 1.  Assume $g=(1,0)$ below the curve $\gamma_1$
and to the right of the curve $\gamma_2$, while $g=(0,1)$ above the
curve $\gamma_1$.  Observe that there exists a Carath\'eodory
solution $x(\cdot)$ of (1.11) going through the points $A$ and $B$.
Next, consider a perturbed solution $x^\ve$, following the vector
field horizontally up to $P$, jumping from $P$ to $P'$,
then moving vertically to $Q$ and horizontally afterwards.
To fix the ideas, assume that 
$$
A=(0,0),\quad B=(0,1),\quad P=(\ve, -\ve^\alpha),\quad
P'=(\ve,\ve^\alpha),
$$
$$
\gamma_1=\big\{ x_2=|x_1|^\alpha\big\},\qquad\qquad
\gamma_2=\big\{ x_1=|x_2-1|^\beta\big\},
$$
In this case the trajectory $x^\ve$ is a solution of a perturbed
system where $\dot w$ is a single Dirac mass of strength
$|P'-P|=2|\ve|^\alpha$.  On the other hand, after both
trajectories have switched to the right of the curve
$\gamma_2$ their distance is
$\|x^\ve-x\|=\ve^{1/\beta}$.  
In this example, the distance between solutions
grows much worse than linearly w.r.t.~the strength of the perturbation,
Indeed, the only estimate available is
$$
\big\|y-x\big\|_{\L^\infty}=\O(1)\cdot 
\big(\tv\{w\}\big)^{1/\alpha\beta}.
\tag 1.13
$$
\midinsert
\vskip 10pt
\centerline{\hbox{\psfig{figure=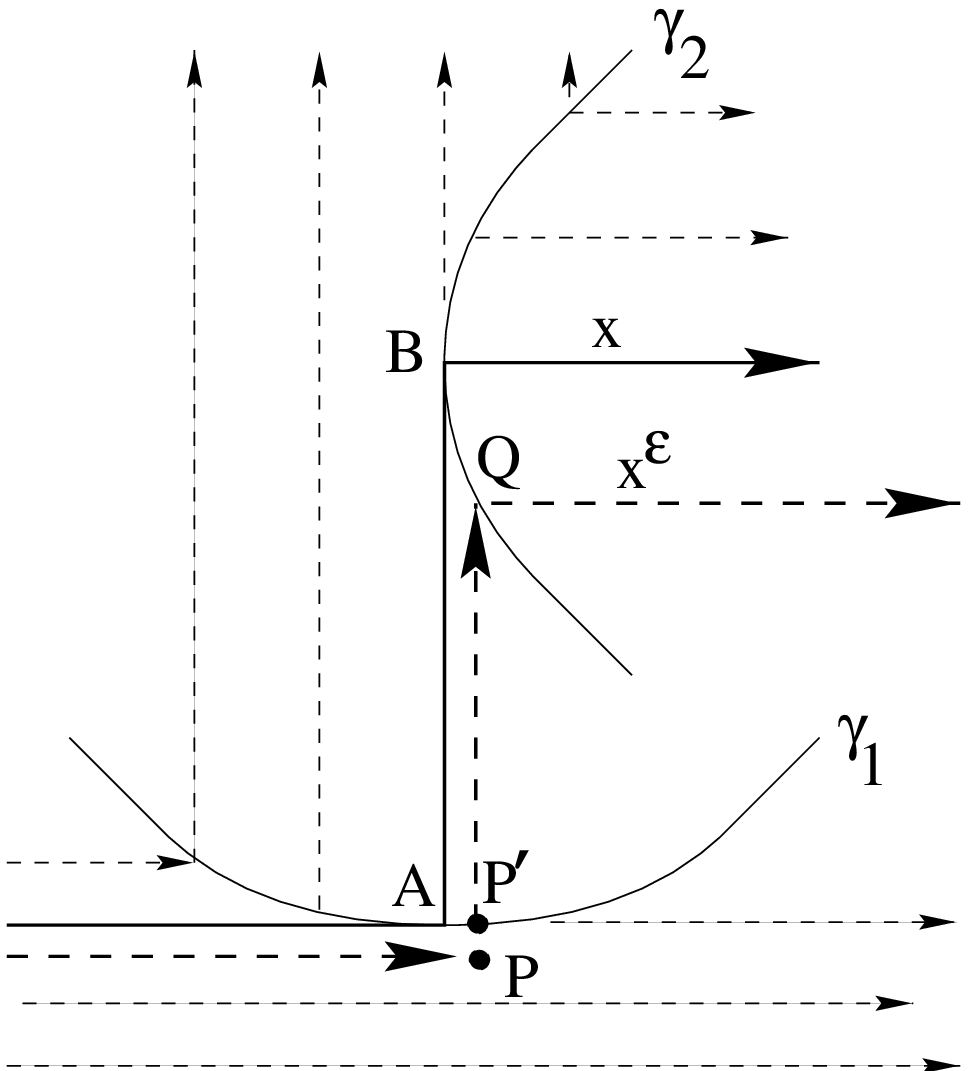,width=6cm}}}
\centerline{\hbox{\smc figure 1}}
\vskip 5pt
\endinsert
One conjectures that the situation is better when the 
patches are in ``generic'' position.  Observe that in (1.13) the 
numbers $\alpha$ and $\beta$ are determined by 
the order of tangency of the curves
$\gamma_1,\gamma_2$ with the vector field $g$.  By an arbitrarily
small displacement of the curves $\gamma_1,\gamma_2$
we can arrange so that there is no trajectory connecting the two 
point of tangency $A$ and $B$ (fig.~2).  
Moreover, we can assume that the
tangency is only of first order.  For generic patchy vector fields on
$\R^2$, in Corollary~1.1 one thus expects an estimate of the form
$$
\big\|y-x\big\|_{\L^\infty}=\O(1)\cdot\big(\tv\{w\}\big)^{1/2}.
$$
\noindent
Here the exponent $1/2$ is due to the fact that first order tangencies
are not removable by small perturbations.  In higher space dimensions, an
even lower exponent is expected.
To obtain an error estimate which is linear w.r.t.~the strength of
the perturbation, one thus needs to remove all these tangencies.
This cannot be achieved if the patches have smooth boundary, but is
quite possible if we allow ``polyhedral'' patches (fig.~3).
\v
\midinsert
\vskip 10pt
\centerline{\hbox{\psfig{figure=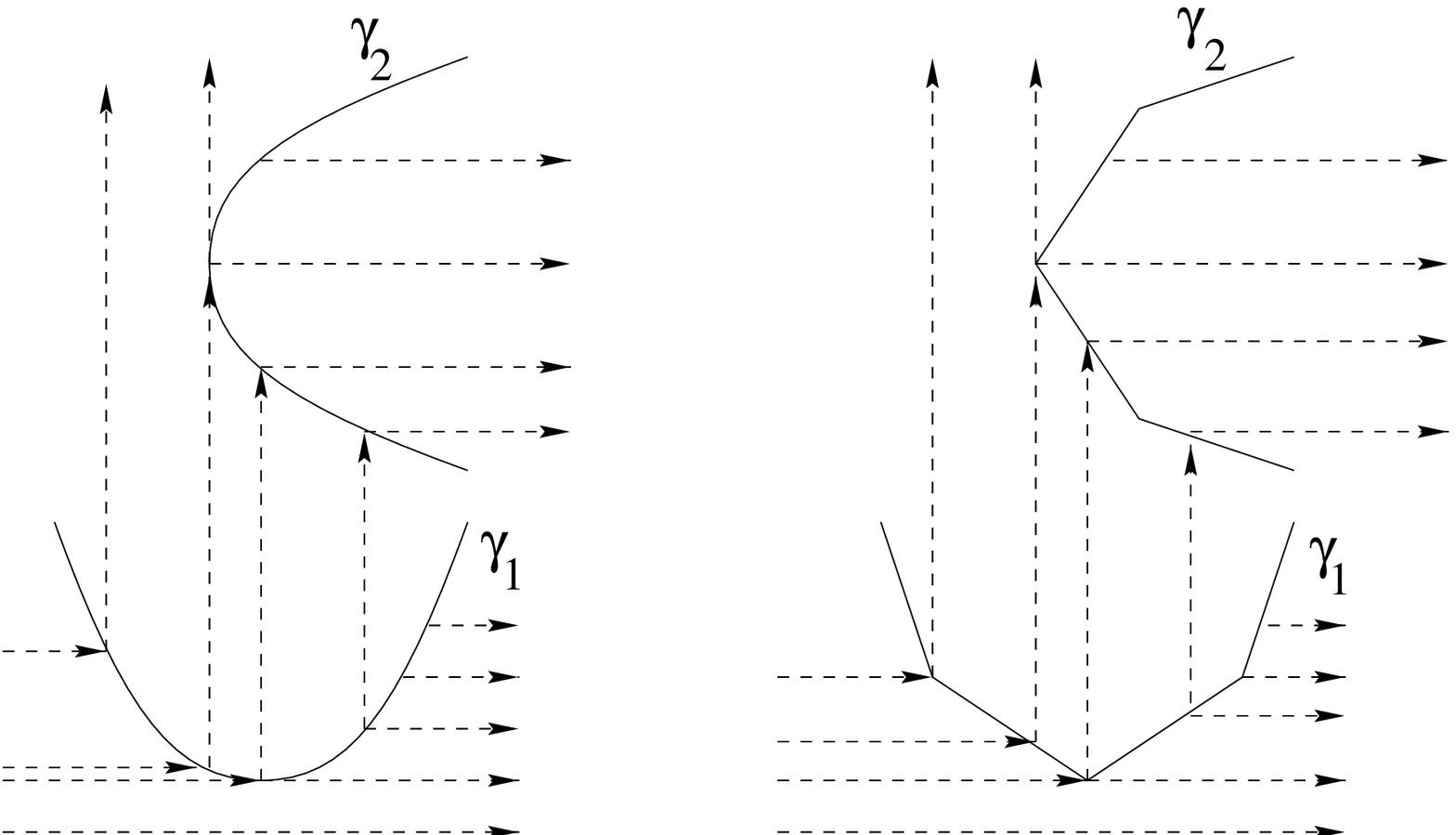,width=10cm}}}
\centerline{\hbox{\smc figure 2}\hskip 1.3in 
\hbox{\smc figure 3}
\hskip 0.6in}
\vskip 10pt
\endinsert

Throughout the following, we write
$d(x,\,A)=\inf\big\{|x-y|~:~ y\in A\big\}$ for 
the distance of a point $x$ from the set $A \subset \Bbb R^n,$
and denote by $\overset \circ \to A$ the 
interior of $A$.
\v
\n{\bf Definition 1.5.} \ Let $\Omega\subset\R^n$ be an open
domain whose boundary is contained in a finite set of
hyperplanes. 
Call
$T_\Om(x)$ the 
tangent  cone to  $\Om$  at the point $x$,
defined by
$$ 
T_\Om(x)\doteq
\bigg\{
v\in \Bbb R^n~:~ \liminf_{t \downarrow 0}
\frac{d\big(x+tv,\ \Om\big)}{t}=0
\bigg\}.
\tag 1.14
$$
We say that a smooth vector field  $g$
defined on a neighborhood of $\ov\Om$ is an {\it inward-pointing
vector field} on \, $\Om$ \, if, 
$$
g(x)\in  \overset \circ \to T_\Om(x)
\qquad\quad \forallt\qquad x\in\partial\Om.
\tag 1.15
$$
\v

\noindent
The pair $\big(\Om,\, g\big)$ will be  called a {\it polyhedral patch}.
\v

Clearly, at any regular point $x\in\partial\Om$, the
interior of the 
tangent cone $T_\Om(x)$ is precisely the set of all vectors
$v\in \Bbb R^n$ that satisfy
$$
\la v,~\bn(x)\ra <0.
$$
and hence (1.15) coincides with the inward-pointing 
condition (1.7).
\v

Replacing ``patches'' with ``polyhedral patches'' 
in Definition~1.2 we obtain the
notion of {\it polyhedral patchy vector field}.  
For such fields,
it is expected that 
impulsive perturbations of the form (1.1) 
should generically produce a
perturbation on the set of trajectories which is of exactly 
the same order of magnitude as the strength of the impulse
on the right hand side.

To avoid lengthy technicalities, we shall consider here 
only 
the planar case, i.e.~{\it polygonal patchy vector fields}.
We conjecture 
that the same result  holds true for generic polyhedral patchy
vector fields on $\Bbb R^n.$

\v
\n{\bf Theorem 1.} {\it For a generic
polygonal patchy vector field $g$ on $\R^2$, 
whose values are bounded away from zero, one has the 
the following stability property.

Given any $T>0$ and any compact set $K\subset \Bbb R^2$ , 
there exist constants
$C$, \, $\delta>0$, 
such that the following holds.
For every solution $y:[0,T]\mapsto\R^2$ of~(1.1)
with $y(0)\in K$ and $\tv\{w\}<\delta$, there exists a solution
$x:[0,T]\mapsto \R^2$ of (1.11) such that}
$$
\big\|x - y\big\|_{\L^\infty([0,T])}\leq C\cdot \tv\{w\}.
\tag 1.16
$$
\v
A precise description of the generic conditions which
guarantee the estimate (1.16) will be
given in Section~2.  
Roughly speaking, one requires that 
the boundary of every patch $\Om_\alpha$ be
transversal to all fields $g_\beta$, with
$\beta\leq\alpha$. 
\v

Throughout the paper, 
by $B(x,r)$ we denote the closed ball centered at
$x$ with radius $r$.  The closure, the interior and the 
boundary of a set
$\Omega$ are written as $\overline\Omega$,
$\overset \circ \to \Omega$ and $\partial\Omega$, respectively.
\v

The paper is organized as follows. 
In Section~2 we introduce a class
of polygonal patchy vector fields for which we
will establish the stability property stated in Theorem~1,
and we show that we can always replace
a solution of the perturbed system (1.1)
with a piecewise smooth concatenation
of solutions of the unperturbed system (1.11),
so that their distance is 
of the same order of magnitude as the impulsive term $\dot w$.
To establish this result contained in Proposition~2.4, we rely 
on two technical lemmas (Lemma~2.2 and Lemma~2.3)
whose rather
lengthy proofs
are postponed to Section~4 (Appendix).
In Section~3 we 
first show in Proposition~3.1 that, 
for every function \, $y(\cdot)$ \, that is a 
concatenation
of two solutions of (1.11) 
(and thus admits a single jump discontinuity), there exists
a solution  \, $x(\cdot)$ \, of (1.11) 
for which the linear estimate (1.16) holds,
and then we complete the proof of Theorem~1
establishing Lemma~3.2.

\vsk

\n{\medbf 2 - Preliminary Stability Estimates.}
\v

Let $\Cal P\Cal P\Cal V\Cal F$ denote the set of all bounded, 
{\bf polygonal patchy
vector fields} 
$\big(g,\  (\Omega_\alpha,~g_\alpha)_{\alpha\in\A}\big)$
on~$\Bbb R^2,$ that are uniformly non-zero.
A {\it condition} \, $P$ \, for a patchy vector field \,
$\big(g,\  (\Omega_\alpha,~g_\alpha)_{\alpha\in\A}\big)\in 
\Cal P\Cal P\Cal V\Cal F$ \,
is a logic proposition that
can be expressed in terms of the fields $g_\alpha$
and (or) the domains~$\Omega_\alpha.$ We write \, $P(g)$ \,
if  \,
$\big(g,\  (\Omega_\alpha,~g_\alpha)_{\alpha\in\A}\big)$ \,
satisfies $P,$ and we say that \, $P$ \, 
is {\it generic} if \, $\{g\in\Cal P\Cal P\Cal V\Cal F~:~P(g)\}$ \,
is a generic subset of $\Cal P\Cal P\Cal V\Cal F.$
\v

We state now a generic condition that yields the linear estimate
(1.16) of the effect of impulsive perturbations on the solutions
of the unperturbed system (1.11). 

\v

\parindent 35pt 
\item{${{\boldkey (}\bold C
{\boldkey )}}$} \ For any given  domain \,
$\Omega_\alpha,$ \, and for any line \, $r_\gamma$ \, containing an edge
of the boundary \, $\partial \Omega_\gamma$ \, of some \, $\Omega_\gamma
, \ \gamma>\alpha,$ \, the field \, $g_\alpha(x)$ \, is transversal to 
\, $r_\gamma$ \, at every point \, $x\in r_\gamma\cap
\overline{\big(\Om_\alpha \setminus \bigcup_{\beta > {\alpha}} 
\Omega_\beta\big)}.$ 
%
\v

\v

\parindent 20pt

\noindent
In this section we will show that, given a 
polygonal patchy vector field $g$ satisfying
condition {\bf (C)},
in order to establish the stability  estimate (1.16)
for an arbitrary solution $t\mapsto y(t)$ of (1.1)
we can always replace $y(\cdot)$ with a piecewise smooth
map $t\mapsto y^\ds(t)$ that is
a concatenation
of solutions of the unperturbed system~(1.11).
This result is contained in Proposition~2.4
and  is based on two technical
Lemmas (Lemma~2.2-2.3) whose proof is postponed to 
Section~4.
Since we shall always consider throughout the
paper solutions of (1.11) or of (1.1) that are contained
in some fixed compact set~$K,$ we will assume
without loss of generality
that every domain \, $\Omega_{\alpha}$ \,
is bounded since, otherwise, one can replace
\, $\Omega_{\alpha}$ \, with its intersection
\, $\Omega_\alpha \cap \Omega'$ \, with 
a polygonal domain 
\, $\Omega'$ \, that contains \, $\overline K$,
preserving the inward-pointing condition (1.15)
and the transversality condition {\bf (C)}.
\v

By the basic properties of a patchy vector field,
for every solution $t\mapsto x(t)$ of (1.11) the corresponding map
$t\mapsto \alpha^*\big(x(t)\big)$ in (1.9) is non-decreasing.
Roughly speaking, a trajectory can move from
a patch $\Omega_\alpha$ to another patch $\Omega_\beta$ only if
$\alpha<\beta$. This property no longer holds in the presence of
an impulsive perturbation. However, it was shown in \ABfs\
that, for a solution $t\mapsto y(t)$ of (1.1), one can slightly modify
the impulsive perturbation $w$, say replacing it by another
perturbation $w^\ds$, so that the map 
$t\mapsto\alpha^*\big(y^\ds(t)\big)$ is monotone
along the corresponding trajectory $t\mapsto y^\ds(t)$.
Namely, the following holds.
\v
\n{\bf Proposition~2.1.} [2, \, Proposition~2.2] \
{\it Let $g$ be a 
patchy
vector field on an open domain~$\Omega \subset \R^n.$ \linebreak
Then, 
given $T>0$ and any compact set $K \subset \Omega,$ \,
there exist
constants $C'=C'(K,T)>0,$ \,  $\delta'=\delta'(K,T)>0,$ 
such that the following holds.

For every BV function \, $w=w(t)$ \, with \, $\tv\{w\}<\delta',$ \,
and for every solution \, $y:[0,T]\mapsto \Omega$  \, of the
Cauchy problem (1.1)-(1.2) with \, $y_0\in K,$ \,
there is a BV function \, $w^\ds=w^\ds(t)$ \,
and a left continuous solution \,
$y^\ds:[0,T]\mapsto \Omega$ \, of 
$$
\dot y^\ds=g(y^\ds)+\dot w^\ds\,,
\tag 2.1
$$
so that the map $t \mapsto \alpha^*(y^\ds(t))$ is 
non-decreasing, and there holds
$$
\aligned
\tv\{w^\ds\}
&\leq C'\cdot \tv\{w\}\,,
\\
\noalign{\medskip}
\big\|y^\ds-y\big\|_{\L^\infty([0,T])}
&\leq C'\cdot \tv\{w\}\,.
\endaligned
\tag 2.2
$$
}
\v

The next  Lemma shows that we can replace 
the solution $t\mapsto y^\ds(t)$ of (2.1) with
a piecewise smooth function  $t\mapsto y^\di(t)$
so that the map 
$t\mapsto\alpha^*\big(y^\di(t)\big)$ is still non-decreasing
and, for every interval
$$
I_\alpha \doteq \big\{t\in [0,T]~;~ \alpha^*(y^\di(t))=
\alpha\big\}\,,
\qquad\ \alpha \in \text{Im}\big(\alpha^*\circ y^\di\big)\,,
$$
\, $y^\di\!\restriction_{I_\alpha}$ \, is a concatenation
of trajectories of (1.11) whose endpoints lie on the edges of the
domain
$$
D_\alpha
\doteq
\Omega_\alpha \setminus 
\displaystyle{\bigcup_{\beta > \alpha} \Omega_\beta}\,.
\tag 2.3
$$

\v
\n{\bf Lemma~2.2.} {\it Let $g$ be a 
uniformly
non-zero polygonal  patchy
vector field on $\R^2,$
associated to a family of polygonal patches
$\big\{ (\Omega_\alpha,~g_\alpha)  ;\, \alpha\in\Cal A\big\}$, and
assume that condition {\bf (C)} is satisfied.
Then, given $T>0$ and any compact set $K \subset \Bbb R^2,$ \,
there exist
constants $C''=C''(K,T), \,  \delta''=\delta''(K,T)>0,$ 
so that, for every BV function $w=w(t)$ with $\tv\{w\}<\delta'',$
and for every solution \, $y: [0,T]\mapsto\Bbb R^2$ \, of~(1.1),
starting at some point \, $y_0\in K,$ \,  
there exists a left continuous, piecewise smooth function \,
$y^\di:[0,T]\mapsto\Bbb R^2$ \, enjoing the properties:}
\v
\itemitem{$a'$)}
{\it The map \, $t\mapsto \alpha^*\big(y^\di(t)\big)$ \, is
non-decreasing.}
\smallskip

\itemitem{$b'$)}
{\it If we let 
$$
\big\{\,
\alpha_{i'_1},\dots,\alpha_{i'_{m^\di}}
\,\big\}
= 
\text{Im}\big(\alpha^* \circ y^\di\big),
\tag 2.4
$$
with
$$
\alpha_{i'_1}< \dots < \alpha_{i'_{m^\di}},
\tag 2.5
$$
and denote 
 \, $D_{\alpha_{i'_k}}$ \, a polygonal
domain defined as in (2.3), then, for every interval
$$
]\tau'_k,\,\tau'_{k+1}]
\doteq
\Big\{
t\in [0,\,T]~:~y^\di(t)\in 
 D_{\alpha_{i'_k}}
\Big\}\,,
$$
there exists a partition 
\, $\tau'_k=t_{k,1}<t_{k,2}<\cdots<t_{k,q_k}=\tau'_{k+1},$ \,
with \, $q_k$ \, less or equal to the number of edges of
the domain \, $D_{\alpha_{i'_k}}$ \, in (2.3),
so that,
on each 
\, $]t_{k,\ell-1},\,t_{k,\ell}\,[,$ \,
the function  \, $y^\di(\cdot)$ \, is a classical solution
of 
%
$$
\dot y= g_{\alpha_{i'_k}}(y)\,.
\tag "$(2.6)_k$"
$$
and the points
\, $y^\di(t_{k,\ell}), \ y^\di(t_{k,\ell}^{\ \ \,+}), \ 
t_{k,\ell} \neq 0,\, T,$ \,
lie on 
different edges of the domain~$D_{\alpha_{i'_k}}.$}
\smallskip

\itemitem{$c'$)} 
\vskip -0.7cm

%
$$
\align
\sum
\Sb
\\
1\leq k\leq m^\di
\\
\noalign{\smallskip}
1\leq \ell< q_k
\endSb
\big|y^\di(t_{k,\ell}^{\ \ \,+})-y^\di(t_{k,\ell})\big|
&\leq C''\cdot \tv\{w\}\,,
\tag 2.7
\\
\noalign{\medskip}
\big\|y^\di-y\big\|_{\L^\infty([0,T])}
&\leq C''\cdot \tv\{w\}\,.
\tag 2.8
\endalign
$$
\v

The next Lemma shows that
for every piecewise smooth function \, $y^\bem(\cdot)$ \,
which is a concatenation
of trajectories of (1.11) and takes values in a domain \, $D_\alpha$ \,
as (2.3), there is a solution \, $x(\cdot)$ \, of 
(1.11) whose $\L^\infty$ distance from $y^\bem(\cdot)$ 
grows linearly with the total amount of jumps in \, $y^\bem.$
\v

\n{\bf Lemma~2.3.} {\it There exist constants \,
$\overline C, \, \overline\delta>0,$ \, so that the
following hold. Let \, $y^\bem :~]\tau_0, \tau_1]\mapsto\Bbb R^2$ \,
be any left-continuous, piecewise smooth
function  
having the properties:
\v
%
%
\itemitem{$a''$)}
the function  \, $y^\bem(\cdot)$ \, is a solution
of \, $\dot y = g_{\alpha_o}(y)$ \, on every interval
\, $]t'_{\ell-1},\,t'_{\ell}\,[\,$ \,
of a partition
\, $t'_1=\tau_0<t'_2<\cdots<t'_{q_o}=\tau_1$ \, of \, $[\tau_0,\tau_1],$ \,
and one has
$$
y^\bem(t)\in D_{\alpha_o}
\qquad\quad\forall~t\in \,]t'_{\ell-1},\,t'_{\ell}\,[
\qquad\quad \forall~k\,,
\tag 2.9
$$
where \, $D_{\alpha_o}$ \, denotes a polygonal
domain defined as in (2.3).
Moreover, one of the following three cases occurs:
\v
\parindent 37pt

\itemitem{1:}
the points
\, $\{y^\bem(t'_{\ell}), \, 1<\ell \leq q_o\}, \ \{y^\bem({t'_{\ell}}^{\, +}), 
\ 1\leq \ell < q_o\},$ \,
lie on 
different edges of the domain~$D_{\alpha_o},$
\smallskip

\itemitem{2:}
the points
\, $\{y^\bem(t'_{\ell}), \, 1<\ell < q_o\}, \ \{y^\bem({t'_{\ell}}^{\, +}), \ 
1\leq \ell < q_o\},$ \,
lie on 
different edges of the domain~$D_{\alpha_o},$
\smallskip

\itemitem{3:}
the points
\, $\{y^\bem(t'_{\ell}), \, 1<\ell \leq q_o\}, \ \{y^\bem({t'_{\ell}}^{\, +}), \ 
1< \ell < q_o\},$ \,
lie on 
different edges of the domain~$D_{\alpha_o};$
\medskip

\parindent 20pt

\itemitem{$b''$)} 
$$
\Delta(y^\bem)\doteq
\sum_{\ell=2}^{q_o-1}\big|y^\bem({t'_\ell}^{\, +})-
y^\bem(t'_\ell)\big|<\overline\delta\,.
\tag 2.10
$$
\v

\noindent
Then, there exist a point \, $Q_{\alpha_o}=Q_{\alpha_o}(y^\bem)
\in \overline D_{\alpha_o},$ \, 
and a time \, $\sigma_{\alpha_o}=\sigma_{\alpha_o}(y^\bem)>0,$ \, 
so that:
\v
\itemitem{$c''$)} 
$$
x^{\alpha_o}\big(t;\ \tau_0,\, Q_{\alpha_o}\big)\in D_{\alpha_o}
\qquad\quad \forall~t\in~]\,\tau_0,\, \sigma_{\alpha_o}[\,,
\tag 2.11
$$
\medskip

\itemitem{$d''$)} if Case~$a''$-{\it 1}) occurs then 
\, $Q_{\alpha_o},\, x^{\alpha_o}\big(\sigma_{\alpha_o};\ \tau_0,\, 
Q_{\alpha_o}\big)\in \partial D_{\alpha_o},$ \,
if Case~$a''$-{\it 2}) occurs then 
\, $Q_{\alpha_o}\in \partial D_{\alpha_o},$ \,
if Case~$a''$-{\it 3}) occurs then 
\, $x^{\alpha_o}\big(\sigma_{\alpha_o};\ \tau_0,\, 
Q_{\alpha_o}\big)\in \partial D_{\alpha_o},$ 
\medskip

\itemitem{$e''$)}
\vskip -0.7cm

$$
\align
\big|\sigma_{\alpha_o}-\tau_1\big|
&\leq \overline C \cdot\Delta(y^\bem)\,,
\tag 2.12
\\
\noalign{\medskip}
\big|x^{\alpha_o}\big(t;\ \tau_0,\, Q_{\alpha_o}\big)- y^\bem(t)\big|
&\leq \overline C \cdot\Delta(y^\bem) \qquad\ \forall~t
\in~]\,\tau_0,\,\min\{\tau_1, \sigma_{\alpha_o}\}]\,.
\tag 2.13
\endalign
$$
}
\v

Relying on Lemmas~2.2-2.3 we are now in the position
to show that, for every solution $t\mapsto y(t)$
of the perturbed system (1.1), we can find a piecewise smooth
map $t\mapsto y^\ds(t)$ that is a
concatenation
of solutions of the unperturbed system (1.11)
and whose $\L^\infty$ distance from $y(\cdot)$ is
of the same order of magnitude as the impulsive term $\dot w$.
\v

\n{\bf Proposition~2.4.} {\it In the same setting of Lemma~2.2, 
given $T>0$ and any compact set $K \subset \Bbb R^2,$ there exist
constants $C'''=C'''(K,T), 
\,  \delta'''=\delta'''(K,T)>0,$ 
so that, for every BV function $w=w(t)$ with $\tv\{w\}<\delta''',$
and for every solution \, $y: [0,T]\mapsto\Bbb R^2$ \, of~(1.1),
starting at some point \, $y_0\in K,$ \, 
there exists a left continuous, piecewise smooth function \,
$y^\ds:[0,T]\mapsto\Bbb R^2$ \, with the 
following properties:
\v

\itemitem{$a'''$)}
The map \, $t\mapsto \alpha^*\big(y^\ds(t)\big)$ \, is
non-decreasing, 
and one has~$y^\ds(0)=y^\ds(0^+).$
\smallskip

\itemitem{$b'''$)}
If we let \, $\alpha_{i'_1}< \cdots < \alpha_{i'_{m^\ds}}$ \, denote the
indices defined for \, $y^\ds(\cdot)$ \, in the same way
as for \, $y^\di(\cdot)$ \, in (2.4)-(2.5), and let
 \, $D_{\alpha_{i'_k}}$ \, denote a 
polygonal
domain defined as in (2.3), then, on every interval
$$
]\tau'_k,\,\tau'_{k+1}]
\doteq
\Big\{
t\in [0,\,T]~:~y^\ds(t)\in 
 D_{\alpha_{i'_k}}
\Big\}\,,
\tag 2.14
$$
the function  \, $y^\ds(\cdot)$ \, is a classical solution
of 
$$
\dot y= g_{\alpha_{i'_k}}(y)\,.
\tag "$(2.15)_k$"
$$
Moreover, one has
$$
y^\ds(\tau_k)\in \partial D_{\alpha_{i'_{k-1}}},
\qquad\quad
y^\ds(\tau_k^{\ +}) \in \partial D_{\alpha_{i'_k}}
\qquad\qquad \forall~1<k\leq m^\ds\,.
\tag 2.16
$$
%
\smallskip

\itemitem{$c'''$)} 
%
$$
\align
\sum_{k=1}^{m^\ds}
\big|y^\ds({\tau'_k}^{\, +})-y^\ds(\tau'_k)\big|
&\leq C'''\cdot \tv\{w\}\,,
\tag 2.17
\\
\noalign{\medskip}
\big\|y^\ds-y\big\|_{\L^\infty([0,T])}
&\leq C'''\cdot \tv\{w\}\,.
\tag 2.18
\endalign
$$
}
\v

\noindent
{\bf Proof.} \ 
Thanks to 
Lemma~2.2, 
in order to establish Proposition~2.4
it will be sufficient 
to show that there exist constants 
\, $C'''=C'''(K), \, 
\delta'''=\delta'''(K)>0,$ \,
so that, for any given piecewise smooth function 
\, $y^\di: [0,T]\mapsto\Bbb R^2,$ \, $y^\di(0)\in K,$
enjoiong properties $a')$, $b')$, $c')$ stated in Lemma~2.2, 
and satisfying the condition
$$
\Delta(y^\di)\doteq
\sum
\Sb
\\
1\leq k\leq m^\di
\\
\noalign{\smallskip}
1\leq \ell< q_k
\endSb
\big|y^\di(t_{k,\ell}^{\ \ \,+})-y^\di(t_{k,\ell})\big|<\delta'''\,,
\tag 2.19
$$
we can construct another 
piecewise smooth function \, $y^\ds: [0,T]\mapsto\Bbb R^2$ \, 
having the properties \,
$a''')$, $b''')$, and satisfying the estimate
$$
\big\|y^\ds-y^\di\big\|_{\L^\infty([0,T])}
\leq C'''\cdot 
\Delta(y^\di)\,.
\tag 2.20
$$
To this purpose, letting \, $\overline C, \ \overline \delta>0$ \,
be constants chosen according with Lemma~2.3,
consider
a piecewise smooth function 
\, $y^\di: [0,T]\mapsto\Bbb R^2,$ \, $y^\di(0)\in K,$
enjoiong properties $a')$, $b')$, $c')$ stated in Lemma~2.2, 
and satisfying~(2.19) with \, 
$$
\delta'''=\overline \delta.
\tag 2.21
$$ 
%
The smoothness
of the vector fields \, $g_\alpha$ \, guarantees that
$y^\di(\cdot)$ \, takes values in some compact
set \, $K'\subset \Bbb R^2$ \, (depending only on $K, \, T, \, 
\overline \delta\,)$. \,
Let \,
$\big\{
\Om_{\alpha_i}~:~i=1,\dots,N
\big\}$ \, be  the collection of polygonal domains
that intersect  \, $K',$ \,
and set
$$
M\doteq 
\sup\big\{|g_{\alpha_i}(y)|~:~y\in \Om_{\alpha_i},
\quad i=1,\dots,N\big\}\,.
\tag 2.22
$$
Let \, $0=\tau'_1<\tau'_2<\cdots<\tau'_{m^\di+1}=T,$ \, 
$m^\di\leq N,$ \,
be the partition of \, $[0, T]$ \, induced by \, $y^\di(\cdot)$ \,
according with property  $b')$, and observe that
every restriction map
\, $y^\di\!\!\restriction_{_{]\tau'_k,\, \tau'_{k+1}]}}\,,$ \,
$1 \leq k \leq m^\di,$ \,
is a piecewise smooth function that enjoys the properties
$a'')$, $b'')$ stated in Lemma~2.3.
Let
$$
\sigma_{\alpha_{i'_k}}\doteq
\sigma_{\alpha_{i'_k}}\big(
y^\di\!\!\restriction_{_{]\tau'_k,\, \tau'_{k+1}]}}\big)\,,
\qquad\quad
Q_{\alpha_{i'_k}}\doteq 
Q_{\alpha_{i'_k}}\big(
y^\di\!\!\restriction_{_{]\tau'_k,\, \tau'_{k+1}]}}\big)\,,
\qquad 1 \leq k \leq m^\di\,,
$$
be the points and times having the properties 
$c'')$, $d'')$, $e'')$ given by Lemma~2.3.
Then, consider the sequence of points 
\, $\tau_1''\doteq
0<\tau_2''<\cdots<\tau_{m''+1}'' \leq T,$ \,
$m''\leq m^\di,$ \, recursively defined by setting
$$
\tau_{k+1}''\doteq
\tau_k''-\tau_k'+
\sigma_{\alpha_{i'_k}}\,,
\tag 2.23
$$
for all \, $1\leq k\leq m^\di$ \, such that \, $\tau_k''-\tau_k'+
\sigma_{\alpha_{i'_k}}<T,$ \, and then letting
$$
\gathered
m''\doteq
\max\big\{
1<k\leq m^\di~:~ \tau_k''-\tau_k'+
\sigma_{\alpha_{i'_k}}<T
\big\}\,,
\\
\noalign{\medskip}
\tau_{m''+1}''\doteq 
\min\big\{T,\
\tau_{m''}''-\tau_{m''}'+
\sigma_{\alpha_{i'_{m''}}}
\big\}\,.
\endgathered
\tag 2.24
$$
Next, letting \, $x^{g}\big(t;\, t_0,\,x_0)$ \,
denote a solution of (1.11) starting 
from \, $x_0$ \, at time \, $t_0,$ \,
define the map
\, $y^\ds: [0,T]\mapsto\Bbb R^2$ \, as follows: \,
$y^\ds(0)\doteq Q_{\alpha_{i'_1}},$ \, and
%
$$
y^\ds(t)\doteq
\cases
x^{\alpha_{i'_k}}\big(t+\tau'_k-\tau''_k;\
\tau'_k,\, Q_{\alpha_{i'_k}}\big)
\quad\quad 
&\forall~t\in ~]\tau''_k,\, \tau''_{k+1}]\,,
\qquad\quad  1\leq k \leq m''\,,
\\
\noalign{\medskip}
x^{g}\big(t;\, \tau''_{m''+1},\,y^\ds(\tau''_{m''+1})\big)
\quad\quad 
&\forall~t\in ~]\tau''_{m''+1},\,T]\,.
\endcases
\tag 2.25
$$
By construction,
the properties $c'')$, $d'')$ of \, 
$\sigma_{\alpha_{i'_k}},\, Q_{\alpha_{i'_k}}$ \,
given by Lemma~2.3,
together with
the general properties of the
solutions of a patchy system 
(recalled in Section~1), guarantee
that  the map \, $t \mapsto y^\ds(t)$ \,
enjoys the properties \,
$a''')$, $b''')$ stated in Proposition~2.4.
Moreover, observe that by property $e'')$ of \,  Lemma~2.3
one has
$$
\aligned
\big|\sigma_{\alpha_{i'_k}}-\tau'_{k+1}\big|
&\leq \overline C \cdot\Delta(y^\di)\,,
\\
\noalign{\medskip}
\big|x^{\alpha_{i'_k}}\big(t;\ \tau'_k,\, Q_{\alpha_{i'_k}}\big)- y^\di(t)\big|
&\leq \overline C \cdot\Delta(y^\di) \qquad\ \forall~t
\in~]\,\tau'_k,\ \min\{\tau'_{k+1}, \sigma_{\alpha_{i'_k}}\}]\,,
\endaligned
\qquad\quad \forall~1\leq k \leq m^\di\,.
\tag 2.26
$$
Thanks to (2.26), and since by definition (2.23)-(2.24) one has
$$
\big|\tau''_{k+1}-\tau'_{k+1}\big|\leq
\big|\tau''_k-\tau'_k\big|+\big|\sigma_{\alpha_{i'_k}}-\tau'_{k+1}\big|
\qquad\quad \forall~1\leq k < m''\,,
$$
proceeding by induction on \, $k\geq 1,$ \, we derive
$$
\big|\tau''_{k+1}-\tau'_{k+1}\big|\leq k\cdot 
\overline C \cdot\Delta(y^\di)
\qquad\quad \forall~1\leq k < m''\,.
\tag 2.27
$$
On the other hand, using  (2.22), (2.26),
and relying on 
property  $b')$ of \, $y^\di(\cdot),$\, we obtain
$$
\align
\big|y^\ds\big({\tau''_{k+1}}^{\!+}\big)-y^\ds(\tau''_{k+1})\big|&=
\big|Q_{\alpha_{i'_{k+1}}}-
x^{\alpha_{i'_k}}\big(\sigma_{\alpha_{i'_k}};\ 
\tau'_k,\, Q_{\alpha_{i'_k}}\big)\big|
\\
\noalign{\smallskip}
&\leq\big|Q_{\alpha_{i'_{k+1}}}-y^\di\big({\tau'_{k+1}}^{\!+}\big)\big|+
\big|y^\di\big({\tau'_{k+1}}^{\!+}\big)-y^\di(\tau'_{k+1})\big|+
\\
&\qquad +\big|y^\di(\tau'_{k+1})-
y^\di\big(\min\{\tau'_{k+1}, \sigma_{\alpha_{i'_k}}\}\big)\big|+
\\
&\qquad +
\big|y^\di\big(\min\{\tau'_{k+1}, \sigma_{\alpha_{i'_k}}\}\big)-
x^{\alpha_{i'_k}}\big(\min\{\tau'_{k+1}, \sigma_{\alpha_{i'_k}}\};\ 
\tau'_k,\, Q_{\alpha_{i'_k}}\big)\big|+
\\
&\qquad +\big|x^{\alpha_{i'_k}}\big(\min\{\tau'_{k+1}, 
\sigma_{\alpha_{i'_k}}\};\ 
\tau'_k,\, Q_{\alpha_{i'_k}}\big)-
x^{\alpha_{i'_k}}\big(\sigma_{\alpha_{i'_k}};\ 
\tau'_k,\, Q_{\alpha_{i'_k}}\big)\big|
\\
\noalign{\smallskip}
&\leq \overline C \cdot\Delta(y^\di)+\Delta(y^\di)+
M\cdot\big|\sigma_{\alpha_{i'_k}}-\tau'_{k+1}\big|+
\overline C \cdot\Delta(y^\di)+
M\cdot\big|\sigma_{\alpha_{i'_k}}-\tau'_{k+1}\big|
\\
\noalign{\smallskip}
&\leq \big(1+2\overline C(1+M)\big)\cdot\Delta(y^\di)\,,
\\
\noalign{\smallskip}
&\hskip 1.3in\forall~1\leq k < m''\,.
\tag 2.28
\endalign
$$
\pagebreak

\n
Hence, thanks to (2.22), (2.26)-(2.28), and by
definition (2.25) of \, $y^\ds(\cdot),$\, we derive
$$
\align
\big|y^\ds(t)-y^\di(t)\big|&\leq 
\big|y^\ds(t)-y^\ds\big(t-\tau'_k+\tau''_k\big)\big|+
\big|y^\ds\big(t-\tau'_k+\tau''_k\big)-y^\di(t)\big|
\\
\noalign{\smallskip}
&\leq \sum_{k=2}^{m''}
\big|y^\ds\big({\tau''_k}^{\,+}\big)-y^\ds(\tau''_k)\big|+
M\cdot \big|\tau''_k- \tau'_k\big|+ 
\big|x^{\alpha_{i'_k}}\big(t;\ \tau'_k,\, Q_{\alpha_{i'_k}}\big)- y^\di(t)\big|
\\
\noalign{\smallskip}
&\leq N\cdot\big(1+3\overline C(1+M)\big)\cdot\Delta(y^\di)\,,
\\
\noalign{\medskip}
&\hskip 0.6in\forall~t
\in~]\,\tau'_k,\ \min\{\tau'_{k+1}, \sigma_{\alpha_{i'_k}}\}]\,,
\qquad 1\leq k \leq m''\,,
\tag 2.29
\endalign
$$
while, in the case \, $\sigma_{\alpha_{i'_k}}<\tau'_{k+1},$ \,
we get
$$
\align
\big|y^\ds(t)-y^\di(t)\big|&\leq 
\big|y^\ds(t)-y^\ds\big(\sigma_{\alpha_{i'_k}}\big)\big|+
\big|y^\ds\big(\sigma_{\alpha_{i'_k}}\big)-
y^\di\big(\sigma_{\alpha_{i'_k}}\big) \big|+
\big|y^\di(t)-y^\di\big(\sigma_{\alpha_{i'_k}}\big)\big|
\\
\noalign{\smallskip}
&\leq \sum_{k=2}^{m''}
\big|y^\ds\big({\tau''_k}^{\,+}\big)-y^\ds(\tau''_k)\big|+
M\cdot \big|\sigma_{\alpha_{i'_k}}-\tau'_{k+1}\big|+
\big|y^\ds\big(\sigma_{\alpha_{i'_k}}\big)-
y^\di\big(\sigma_{\alpha_{i'_k}}\big) \big|+
\\
&\qquad\quad+
\Delta(y^\di)+
M\cdot\big|\sigma_{\alpha_{i'_k}}-\tau'_{k+1}\big|
\\
\noalign{\smallskip}
&\leq 3N\cdot\big(1+3\overline C(1+M)\big)\cdot\Delta(y^\di)\,,
\\
\noalign{\medskip}
&\hskip 0.6in\forall~t
\in~]\,\sigma_{\alpha_{i'_k}},\ \tau'_{k+1}]\,,
\qquad 1\leq k \leq m''\,.
\tag 2.30
\endalign
$$
Thus, (2.29)-(2.30) together, yield
$$
\big|y^\ds(t)-y^\di(t)\big| \leq 
3N\cdot\big(1+3\overline C(1+M)\big)\cdot\Delta(y^\di)
\qquad\quad \forall~t
\in~[\,0,\,\tau'_{m''+1}]\,.
\tag 2.31
$$
On \, the \, other \, hand, \, in \, the \, case 
\, $\tau'_{m''+1}<T,$ \,
by \, definition \, (2.24) \, one \, has \, $m''=m^\di,$ \linebreak
$T<\tau''_{m''}-\tau'_{m''}+\sigma_{\alpha_{i'_{m''}}},$ \,
and hence, using (2.26)-(2.27), we get
$$
\align
\big|T-\tau'_{m''+1}\big|
&\leq\big|\tau''_{m''}-\tau'_{m''}+\sigma_{\alpha_{i'_{m''}}}
-\tau'_{m''+1}\big|
\\
&\leq \big|\tau''_{m''}-\tau'_{m''}\big|+
\big|\sigma_{\alpha_{i'_{m''}}}
-\tau'_{m''+1}\big|
\\
&\leq (N+1)\cdot\overline C\cdot\Delta(y^\di)\,.
\tag 2.32
\endalign
$$
Therefore, from (2.28), (2.31)-(2.32) we derive
$$
\align
\big|y^\ds(t)-y^\di(t)\big| &\leq 
\big|y^\ds(t)-y^\ds(\tau'_{m''+1})\big|+
\big|y^\ds(\tau'_{m''+1})-y^\di(\tau'_{m''+1})\big|+
\big|y^\di(t)-y^\di(\tau'_{m''+1})\big|
\\
\noalign{\smallskip}
&\leq\sum_{k=2}^{m''}
\big|y^\ds\big({\tau''_k}^{\,+}\big)-y^\ds(\tau''_k)\big|+
M\cdot \big|T-\tau'_{m''+1}\big|+
\\
&\qquad+\big|y^\ds(\tau'_{m''+1})-y^\di(\tau'_{m''+1})\big|+
\Delta(y^\di)+M\cdot \big|T-\tau'_{m''+1}\big|
\\
\noalign{\smallskip}
&\leq 4N\cdot\big(1+4\overline C(1+M)\big)\cdot\Delta(y^\di)
\\
\noalign{\medskip}
&\hskip 0.6in\forall~t\in~[\tau'_{m''+1},\,T]\,.
\tag 2.33
\endalign
$$
\v
\noindent
Hence, (2.32)-(2.33) together, show that \, $y^\ds(\cdot)$ \,
satisfies the estimates (2.20) 
taking 
$$
C'''> 4N\cdot\big(1+4\overline C(1+M)\big)\,.
\tag 2.34
$$
which completes the proof of  Proposition~2.4.
\fine

\vsk
\n{\medbf 3 - Proof of Theorem~1.}
\v

In view of Proposition~2.4, it is useful to introduce the following
\v

\n{\bf Definition 3.1.} \ A left-continuous,
piecewise smooth function \,
$y^\ds: [0,T]\mapsto\Bbb R^2$ \,  that 
enjoys the properties   $a''')$-$b''')$
stated 
in Proposition~2.4, 
is called
a {\it concatenation of classical solutions \linebreak
(CCS) of}~(1.11).
\v

Notice that, in particular, any Carath\'eodory solution of (1.11)
is always a CCS. 
Before giving the complete proof of Theorem~1
we will first show that, for every given CCS of (1.11) \, $y(\cdot)$ \,
admitting a single jump discontinuity, there exists
a (Carath\'eodory) solution \, $x(\cdot)$ \, of (1.11) 
for which the linear estimate (1.16) holds. Namely, we shall
prove
\v
\n{\bf Proposition~3.1.} {\it Let $g$ be a 
uniformly
non-zero polygonal  patchy
vector field on $\R^2,$
associated to a family of polygonal patches
$\big\{ (\Omega_\alpha,~g_\alpha)  ~; ~\alpha\in\Cal A\big\}$,
and assume that condition {\bf (C)} is satisfied.
Then, given $T>0$
and any compact set $K \subset \Bbb R^2,$ there exist constants \,
$C^{\text{\sl \i v}}=C^{\text{\sl \i v}}(K,T),$ \linebreak 
$\delta^{\text{\sl \i v}}=
\delta^{\text{\sl \i v}}(K,T)>0$ \,
so that the following hold.

\n
Let \, $y^\bem :~[\,\tau_0, \tau_1]\mapsto K,$ \,
\, $y^\beq :~]\,\tau_1, \tau_2]\subset ]\,\tau_1, \tau_0+T]
\mapsto K,$ \, 
be two continuous maps
having the properties:
\v
\itemitem{$i$)}
the function  \, $y^\bem(\cdot)$ \, is a Carath\'eodory solution of (1.11)
and, letting
$$
\alpha_1\doteq\max\big\{
\alpha~;~ \alpha\in \text{Im}(\alpha^* \circ y^\bem)
\big\}\,,
$$
one has
$$
y^\bem(\tau_1)\in \partial D_{\alpha_1}
\tag 3.1
$$
(where \, $D_{\alpha_1}$ \, denotes a polygonal
domain defined as in (2.3));
\smallskip
\itemitem{$ii$)} 
the function  \, $y^\beq(\cdot)$ \, is a solution
of 
\, $\dot y = g_{\alpha_2}(y),$ \, 
for some \, $\alpha_2>\alpha_1,$ \,
and one has
$$
\gathered
y^\beq(t)\in D_{\alpha_2} 
\qquad \ \forall~t\in~]\tau_1, \tau_2]\,,
\\
\noalign{\smallskip}
y^\beq(\tau_1^{\,+})\in \partial D_{\alpha_2}\,;
\endgathered
\tag 3.2
$$
\smallskip
\itemitem{$iii$)} 
$$
\big|y^\bem(\tau_1) -y^\beq(\tau_1^{\,+})\big|< \delta^{\text{\i v}}.
\tag 3.3
$$
\v

\noindent
Then, there exists a Carath\'eodory solution of (1.11) \,
$\Phi_{_{\bem,\beq}}\!\doteq\Phi\big[y^\bem,y^\beq\big]:~
[\,\tau_0,\, \sigma_{_{\bem,\beq}}] \longmapsto \Bbb R^2,$ \,
such that
$$
\align
\big|y^\bem(t)-\Phi_{_{\bem,\beq}}(t)\big|
&\leq C^{\text{\i v}}\cdot
\big|y^\bem(\tau_1) -y^\beq(\tau_1^{\,+})\big|
\qquad\quad\forall~t\in~[\,\tau_0,\,
\min\{\tau_1,\, \sigma_{_{\bem,\beq}}\}]\,,
\tag 3.4
\\
\noalign{\medskip}
\big|y^\beq(t)-\Phi_{_{\bem,\beq}}(t)\big|
&\leq C^{\text{\i v}}\cdot
\big|y^\bem(\tau_1) -y^\beq(\tau_1^{\,+})\big|
\qquad\quad\forall~t\in~]\,\tau_1,\,
\min\{\tau_2,\, \sigma_{_{\bem,\beq}}\}]\,,
\qquad\text{if}\quad \
\sigma_{_{\bem,\beq}}>\tau_1.
\tag 3.5
\endalign
$$
%
%
%
%
%
Moreover, one has
$$
\gathered
\alpha^*\big(\Phi_{_{\bem,\beq}}(t)\big)\leq \alpha_2
\qquad\quad\forall~t\in [\,\tau_0,\, \sigma_{_{\bem,\beq}}],
\\
\noalign{\medskip}
y^\beq(\tau_2)\in \partial D_{\alpha_2}
\qquad\Longrightarrow\qquad
\Phi_{_{\bem,\beq}}(\sigma_{_{\bem,\beq}})
\in \partial D_{\alpha_2}\,,
\endgathered
\tag 3.6
$$
and there holds}
$$
\big|\sigma_{_{\bem,\beq}}-\tau_2\big|
\leq C^{\text{\i v}}\cdot
\big|y^\bem(\tau_1) -y^\beq(\tau_1^{\,+})\big|\,.
\tag 3.7
$$
%
\v

\noindent
{\bf Proof.}  \ 

\noindent 
{\smbf 1.} \ Fix $T>0$ and a compact set \, $K\subset \Bbb R^2.$ \ 
Every trajectory $y : [t_0,\, t_1]\subset
[t_0,\, t_0+T]\mapsto\Bbb R^2$ \,  of~(1.11),
starting at a  point \, $y_0\in K,$ \,
takes values in some compact 
set \, $K'.$ \,
Consider the set of indices
$$
\Cal A_{K'}\doteq \big\{\alpha \in \Cal A~~;~\Omega_\alpha \cap K'
\neq \emptyset
\big\}\,.
\tag 3.8
$$
Let \, $N=|\Cal A_{K'}|$ \, be the number of elements in \, $\Cal A_{K'},$ \,
and set
$$
M\doteq 
\sup\big\{|g_\alpha(y)|~:~y\in \Om_\alpha,
\qquad \alpha\in \Cal A_{K'} \big\}\,.
\tag 3.9
$$
For each \, $\alpha\in \Cal A_{K'},$ \, denote by \, $V_\alpha$ \,
the set of vertices of the polygonal domain \, 
$D_\alpha.$ \, 
Notice that, since \, $g$ \, is a patchy vector field satisfying
condition {\bf (C)},  the Cauchy problem
(1.11)-(1.2) has a unique local forward (Carat\'eodory) solution in the case
\, $y_0\in \big(\bigcup_{\alpha\in\Cal A_{K'}} \Omega_\alpha\big)
\setminus \big(\bigcup_{\alpha\in\Cal A_{K'}} V_\alpha\big),$ \, 
and at most $N$ local forward solutions if 
\, $y_0\in  \bigcup_{\alpha\in\Cal A_{K'}} V_\alpha.$ \,
On the other hand, by the properties 
of the solutions of a patchy system recalled in Section~1, 
the Cauchy problem for (1.11)
has always backward uniqueness.
Therefore, the set~$\Cal T_\alpha$ \, of all graphs of maximal 
(Carat\'eodory) trajectories
of~(1.11) that 
go through some vertex in~$\cup_{\alpha\in \Cal A_{K'}} V_\alpha,$ \, 
and are contained in \, 
$\cup_{\alpha\in \Cal A_{K'}} \overline\Omega_\alpha,$ \,
is finite.
For convenience, with a slight abuse of notation,
we will often write \, $\gamma \in \Cal T_\alpha$ \,
to mean \, $\text{Im}(\gamma)\in \Cal T_\alpha$. \,
The unique backward solution of the
Cauchy problem (1.11)-(1.2), whenever does exist, will be denoted
by \, $t \mapsto x^g(t;\, t_0,\, y_0),$ \, $t\leq t_0.$\ \linebreak
We assume that every vector field $g_\alpha$ is
defined on a neighborhood of \, $\overline\Omega_\alpha$ \,
and we denote, as usual, by 
\, $t \mapsto x^\alpha\big(t;\, t_0,\, x_0\big)$ \
the
solution of the Cauchy 
problem \, 
$$
\dot x=g_\alpha(x), \qquad x(t_0)=x_0\in \overline\Omega_\alpha.
\tag 3.10
$$ \,
By well-posedness of (3.10), there will be some constant  $c_0>1$
so that 
$$
\qquad\quad
\big|
x^\alpha(t;\, t_0,\, x_0)-
x^\alpha(t';\, t_1,\, x_1)
\big|\leq c_0
\big\{|t-t'|+|t_0-t_1|+|x_0-x_1|\big\}
\qquad \quad \forall~x_0,\,x_1\in \overline\Omega_\alpha.
\tag 3.11
$$
%
For every \, $x_0\in 
\partial D_\alpha,$ \, 
we let \, $t_\alpha^+(x_0), \, t_\alpha^-(x_0)$ \,
denote
the time that is necessary to
reach the set \, 
$\overline \Omega_\alpha\setminus D_\alpha$  \, 
starting from $x_0$ 
and following, respectively,  the forward
and backward flow of the vector field $g_\alpha,$ \, i.e. 
$$
\aligned
t_\alpha^+(x_0)
& \doteq
\inf\big\{t>0~;~ x^\alpha(t;\, 0,\, x_0)
\in \overline \Omega_\alpha \setminus D_\alpha
\big\}\,,
\\
\noalign{\smallskip}
t_\alpha^-(x_0)
& \doteq
\sup\big\{t<0~;~ x^\alpha(t;\, 0,\, x_0)
\in \overline \Omega_\alpha \setminus D_\alpha
\big\}\,.
\endaligned
\tag 3.12
$$
%
Using the quantities in (3.12) we define the sets 
of {\it incoming} and {\it outgoing boundary points}
$$
\aligned
\partial_{_{\Cal I}} D_\alpha
&\doteq
\big\{x\in\partial D_\alpha~;~ t_\alpha^+(x)>0
\big\}\,,
\\
\noalign{\smallskip}
\partial_{_{\Cal O}} D_\alpha
&\doteq
\big\{x\in\partial D_\alpha~;~ t_\alpha^-(x)<0
\big\}\,,
\endaligned
\tag 3.13
$$
%
that clearly consist of all the points in
the boundary \, $\partial D_\alpha$ \, where the field $g_\alpha$ is pointing, respectively,
towards the
interior and towards the exterior of \, $D_\alpha.$ \,
%
%
Moreover, for any pair of indices \, 
$\alpha,\beta\in\Cal A_{K'},$ \, $\alpha<\beta,$ \, 
define the set
$$
G_{\alpha,\,\beta}\doteq
\partial_{_{\Cal O}} D_\alpha\cap\partial_{_{\Cal I}} D_\beta\cap
\Big(V_\alpha\cup V_\beta \cup 
\bigcup_{\gamma \in \Cal T_\beta}
\text{Im}(\gamma)\Big)\,.
%
\tag 3.14
$$
Since \, $g_{\alpha}$ \, are smooth, uniformly non-zero 
vector fields
that satisfy the  inward-pointing
condition (1.17) and
the transversality condition {\bf (C)}, and by the properties
of the solutions of a patchy system recalled in Section~1,
one can easily verify
that the following properties hold.
\v

\parindent 25pt

\item{{\bf P1)}} \  {\it There exist  constants \, $c_1,\,
\delta_1>0$ \
(depending only on $K$) so that, given any
domain \, $D_\alpha,$ \, $\alpha\in\Cal A_{K'},$ \, one has:
\vvs
\itemitem{{\it a)}} \ For any \, $x,\, y\in \partial_{_{\Cal I}} D_\alpha$ \,
belonging to the same connected component of \, $\partial D_\alpha\setminus
\bigcup_{\gamma \in \Cal T_\alpha}
\text{Im}(\gamma),$ \, there holds
$$
|x-y| < \delta_1 \qquad\Longrightarrow \qquad \big|t_\alpha^+(x)-t_\alpha^+(y)\big| 
< c_1 \cdot |x-y|.
\tag 3.15
$$
\vvs
\itemitem{{\it b)}} \ Given any \, $x_0\in \partial_{_{\Cal I}} D_\alpha \cap 
\bigcup_{\gamma \in \Cal T_\alpha}
\text{Im}(\gamma),$ \, and any \, $x\in \partial_{_{\Cal I}} D_\alpha \cap B(x_0,\, \delta_1),$ \, 
$x\neq x_0,$ \, there exists \, $\widehat\tau_{x_0,x}^\alpha\in[0,\, t_\alpha^+(x_0)]$ \, such that
%
$$
x^{\alpha}\big(\widehat\tau_{x_0,x}^\alpha;\, 0,\, x_0\big)\in \partial D_\alpha,
\qquad\quad
\big|t_\alpha^+(x)-\widehat\tau_{x_0,x}^\alpha\big|
< c_1 \cdot |x-x_0|.
\tag 3.16
$$
}
\v

\parindent 25pt

\item{{\bf P2)}} \ {\it There exist  constants \, $c_2,\,
\delta_2>0$ \
(depending only on $K$) so that, given any pair of
domains \, $D_\alpha,\, D_\beta,$ \ $\alpha<\beta,$ \, 
$\alpha,\beta\in\Cal A_{K'},$ \,
with \, $\partial_{_{\Cal O}} D_\alpha\cap\partial_{_{\Cal I}} 
D_\beta \neq \emptyset,$ \, 
%
%
the following holds.
For every \, $x_0\in  G_{\alpha,\,\beta},$ \, 
and for any backward (Carat\'eodory) trajectory \,
$t \mapsto x^g(t;\, 0,\, x),$ \ $t\in [-\widehat\tau_x,\,0],$ \, 
of~(1.11) \linebreak
arriving in some point \, $x\in  \partial_{_{\Cal O}} 
D_\alpha \cap B(x_0,\, \delta_2),$ \,
starting from \, $x^g(-\widehat\tau_x;\, 0,\, x)\in 
\bigcup_{\gamma\leq \alpha} \partial D_\gamma,$ \,
and contained in \, $K,$ \,
%
%
%
there is another  backward trajectory
$t \mapsto x^g(t;\, 0,\, x_0),$ \ $t\in [-\widehat\tau_{x_0,x}^g,\,0],$ \, of (1.5)
arriving in \, $x_0,$ \, and such that
$$
\allowdisplaybreaks
\gather
x^g\big(-\widehat\tau_{x_0,x}^g;\, 0,\, x_0\big)\in \bigcup_{\gamma\leq \alpha} \partial D_\gamma,
\tag 3.17
\\
\noalign{\medskip}
\big|\widehat\tau_{x_0,x}^g-\widehat\tau_x\big|< c_2 \cdot |x-x_0|,
\tag 3.18
\\
\noalign{\bigskip}
\quad
\big| x^g\big(t-\widehat\tau_{x_0,x}^g;\, 0,\, x_0\big)- x^g\big(t-\widehat\tau_x;\, 0,\, x\big)\big|
< c_2 \cdot |x-x_0|
\qquad\forall~t\in[0,\ \min\{\widehat\tau_{x_0,x}^g,\, \widehat\tau_x\}].
\tag 3.19
\endgather
$$
}
\parindent 20pt
\v

\noindent
Next, choose \, $\overline \lambda>0$ \, so that,
for any pair of indices \, $\alpha,\beta\in\Cal A_{K'},$ \, $\alpha<\beta,$ \, 
one has
$$
B(x_0,\,\overline\lambda\,) \cap B(y_0,\,\overline\lambda\,)= \emptyset
\qquad\quad\forall~x_0,\,y_0\in G_{\alpha,\beta},\qquad x_0\neq y_0.
\tag 3.20
$$
%
Then,  for any \, 
$0\leq \lambda \leq \overline \lambda,$ \,
let \, $R^1_{\alpha,\beta}(\lambda), \dots, R^{r_{\!_{\alpha,\beta}}}_{\alpha,\beta}(\lambda)$ \,
denote the connected components of \, 
$$
\big(\partial D_\alpha \cup \partial D_\beta\big)\setminus 
\bigcup_{P\in G_{\alpha,\,\beta}} B(P,\, \lambda),
$$ \,
and
set
$$
\rho_{\alpha,\beta}(\lambda)\doteq 
\cases
\min
\Big\{ d\big(R^s_{\alpha,\beta}(\lambda),\,
R^\ell_{\alpha,\beta}(\lambda)\big)
\ \ ;\ \, 1 \leq s,\, \ell \leq r_{\alpha,\beta},\ \ \ s\neq \ell
\Big\}
&\quad\ \text{if} \qquad
\partial_{_{\Cal O}} D_\alpha\cap\partial_{_{\Cal I}} D_\beta
\neq \emptyset\,,
\\
\noalign{\medskip}
d\big(\partial_{_{\Cal O}} D_\alpha,\ \partial_{_{\Cal I}} D_\beta\big)
&\quad\ \text{otherwise.}
\endcases
\tag 3.21
$$
Since by construction  one has
$$
\qquad
\inf
\left\{\,
\frac{\rho_{\alpha,\beta}
(\lambda)}{\lambda}~:~ 0<\lambda \leq \overline \lambda
\,\right\}>0
\qquad  \forall~\alpha,\beta\in \Cal A_{K'},
\quad \alpha<\beta,
\qquad\text{s.t.}\quad 
\partial_{_{\Cal O}} D_\alpha\cap\partial_{_{\Cal I}} D_\beta
\neq \emptyset,
\tag 3.22
$$
there will be
constants \, $c_3>1, \ 0<\delta_3<\overline \lambda/(2c_3), $ \,
so that
$$
\rho_{\alpha,\beta}
(c_3 \cdot \delta) > 2 \delta
\qquad\quad \forall~0<\delta\leq \delta_3.
\qquad  \forall~\alpha,\beta\in \Cal A_{K'},
\quad \alpha<\beta.
\tag 3.23
$$
%
\vs

\noindent
{\smbf 2.} \ Consider now two continuous maps
\, $y^\bem :~]\tau_0, \tau_1]\mapsto K,$ \,
\, $y^\beq :~]\tau_1, \tau_2]\mapsto K,$ \,
having the proper-\linebreak
ties~{\it i)-iii)} with
$$
\delta^{\text{\sl \i v}} \leq \min\left\{\frac{\delta_1}{2 c_3},\, 
\frac{\delta_2}{2 c_3},\ \delta_3
\right\}.
\tag 3.24
$$
To fix the ideas, we shall assume also that
$$
y^\bem(\tau_0^{\,+})\in \bigcup_{\gamma \leq \alpha_1}\partial D_\gamma,
\qquad\quad
y^\beq(\tau_2)\in \partial D_{\alpha_2}.
\tag 3.25
$$
The cases where \, $y^\bem(\tau_0^{\,+})\in 
\bigcup_{\gamma \leq \alpha_1}\overset 
\,\circ \to D_\gamma,$ \, 
or \, $y^\beq(\tau_2)\in\overset 
\,\circ \to D_{\alpha_2},$ \, can be treated in entirely
similar manner.
Set
$$
x^\bem\doteq y^\bem(\tau_1),
\qquad\ \
x^\beq\doteq y^\beq(\tau_1^{\,+}),
\qquad\qquad \Delta\doteq \big|x^\bem-x^\beq\big|,
\tag 3.26
$$ \,
and observe that, since the properties $i)$-$ii)$ of \, $y^\bem(\cdot),\, y^\beq(\cdot)$ \,
imply 
$$
x^\bem \in \partial_{_{\Cal O}} D_{\alpha_1},
\qquad\quad
x^\beq \in \partial_{_{\Cal I}} D_{\alpha_2},
$$
by the definition (3.21) of \,
$\rho_{\alpha_1,\alpha_2}$ \, and because of (3.3), (3.23), (3.24), we deduce
that, if \linebreak
\, $\partial_{_{\Cal O}} D_{\alpha_1}\cap\partial_{_{\Cal I}} D_{\alpha_2}
=\emptyset,$ \, then
$$
\align
\Delta&\geq d\big(\partial_{_{\Cal O}} D_{\alpha_1},\ 
\partial_{_{\Cal I}} D_{\alpha_2}
\big)
= \rho_{\alpha_1,\alpha_2}(c_3\cdot \Delta)
\\
\noalign{\smallskip}
&
> 2 \Delta
\endalign
$$
which yields a contradiction. 
Therefore, it must be \, $\partial_{_{\Cal O}} D_{\alpha_1}\cap\partial_{_{\Cal I}} D_{\alpha_2}
\neq\emptyset$ \,
which, in turn, by definition~(3.14) 
implies \, $G_{\!_{\alpha_1,\alpha_2}} \neq \emptyset.$ \,
In order to construct the Carath\'eodory solution of (1.11) 
\, $\Phi_{_{\bem,\beq}}$ \,
satisfying (3.4)-(3.7),
we will handle separately the case in which the endpoints 
\, $x^\bem, \, x^\beq$ \,
lie on the same connected
component of \, 
$$
\big(\partial D_{\alpha_1} \cup \partial D_{\alpha_2}\big)\setminus  \!\!
\bigcup_{x_0\in G_{\!_{\alpha_1,\alpha_2}}} \!\!\! \!B(x_0,\, c_3\cdot \Delta),
\tag 3.27
$$
and the case
where \, $x^\bem,\,x^\beq$ \, 
belong to the ball \, $B(x_0,\, c_3\cdot \Delta)$ \,
centered at some point \, $x_0\in G_{\alpha_1, \alpha_2}.$ \,
\vs
\pagebreak

\noindent
{\smbf 3.} \ 
{\smc Case~1.} \ Assume that
$$
\aligned
x^\bem
&\in \partial D_{\alpha_1} \setminus\!\!
\bigcup_{x_0\in G_{\!_{\alpha_1,\alpha_2}}}\!\!\! \! 
B(x_0,\, c_3\cdot \Delta),
\\
\noalign{\medskip}
x^\beq
&\in \partial D_{\alpha_2} \setminus\!\!
\bigcup_{x_0\in G_{\!_{\alpha_1,\alpha_2}}}\!\!\! \! 
B(x_0,\, c_3\cdot \Delta),
\endaligned
\tag 3.28
$$
and let \, $R^s_{\alpha_1,\alpha_2}(c_3\cdot \Delta),\,
R^\ell_{\alpha_1,\alpha_2}(c_3\cdot \Delta)$ \, 
be the connected components of 
the set in (3.27) that contain, respectively, \, $x^\bem$ \,  
and \, $x^\beq.$ \, 
Observe that, if \, $s\neq \ell,$ \, then, by the definition (3.21) of \,
$\rho_{\alpha_1,\alpha_2}$ \, and because of (3.3), (3.23), (3.24), we deduce
$$
\align
\Delta&\geq d\big(R^s_{\alpha_1,\alpha_2}(c_3\cdot \Delta),\ 
R^\ell_{\alpha_1,\alpha_2}(c_3\cdot \Delta)
\big)
\\
\noalign{\smallskip}
&\geq \rho_{\alpha_1,\alpha_2}(c_3\cdot \Delta)
\\
\noalign{\smallskip}
&> 2 \Delta
\endalign
$$
which yields a contradiction. Therefore it must be \, $s=\ell,$ \,
i.e. \, $x^\bem, \,  x^\beq$ \, 
lie on the same connected
component of the set in (3.27) and hence one has
$$
x^\bem,\, x^\beq\in
R^s_{\alpha_1,\alpha_2}(c_3\cdot \Delta)
\subset \partial_{_{\Cal O}} D_{\alpha_1}\cap\partial_{_{\Cal I}} D_{\alpha_2}\,.
$$
But then, since (3.3), (3.24) together imply \, $\Delta<\delta_1,$ \,
applying property {\bf P1}-{\it a)} we derive
$$
\big|t_{\alpha_2}^+\big(x^\bem
\big)-t_{\alpha_2}^+\big(x^\beq
\big)\big| 
< c_1 \cdot \Delta\,.
\tag 3.29
$$
On the other hand, since 
\, $y^\beq(\cdot)$ \, satisfies property $ii),$   and because (3.28)
implies  
$$
x^\beq
\in \partial_{_{\Cal I}} D_{\alpha_2}\setminus
\bigcup_{\gamma \in \Cal T_{\alpha_2}}\text{Im}(\gamma)\,,
\tag 3.30
$$
using also $(3.25)_2$ we deduce
$$
\gather
t_{\alpha_2}^+\big(x^\beq
\big)
=\tau_2-\tau_1,
\tag 3.31
\\
\noalign{\medskip}
y^\beq(t)
=x^{\alpha_2}\big(t;\, \tau_1,\,x^\beq
\big)
\qquad\quad\forall~t\in~]\tau_1,\,\tau_2]\,.
\tag 3.32
\endgather
$$
Thus, setting
$$
\sigma_{_{\bem,\beq}}\doteq \tau_1+t_{\alpha_2}^+\big(x^\bem
\big),
\tag 3.33
$$
from (3.29), (3.31) we derive
$$
\gathered
x^{\alpha_2}\big(\sigma_{_{\bem,\beq}};\, \tau_1,\,x^\bem\big)
\in \partial D_{\alpha_2}\,,
\\
\noalign{\medskip}
\big|\sigma_{_{\bem,\beq}}-\tau_2\big|
<c_1\cdot\Delta\,,
\endgathered
\tag 3.34
$$
and, thanks to (3.11), we obtain
$$
\big|y^\beq(t)-x^{\alpha_2}\big(t;\, \tau_1,\,x^\bem
\big) \big|
< c_0 \cdot \Delta \qquad\quad\forall~t\in~]\tau_1, 
\min\{\tau_2, \sigma_{_{\bem,\beq}}\}]\,.
\tag 3.35
$$
Then, define the map 
$$
\Phi_{_{\bem,\beq}}(t)\doteq
\cases
y^\bem(t) \quad &\forall~ t\in ~[\,\tau_0,\,\tau_1]\,,
\\
\noalign{\medskip}
x^{\alpha_2}\big(t;\, \tau_1,\,x^\bem
\big)
\quad &\forall~
t\in ~]\,\tau_1, \sigma_{_{\bem,\beq}}]\,,
\endcases
\tag 3.36
$$
and observe that, by construction, \, $\Phi_{_{\bem,\beq}}(\cdot)$ \, 
is a solution of (1.11)
verifying (3.6). Moreover, (3.4)~trivially holds, while
from (3.34)-(3.35) we recover the
estimates (3.5), (3.7) taking the constant 
\, $C^{\text{\i v}}>\max\{c_0,\,c_1\}.$ \,
Thus,
(3.36) provides the map required by the Proposition
whenever (3.28) is verified.
%

\vs

\noindent
{\smbf 4.} \ 
{\smc Case~2.} \ Assume that 
$$
x^\bem
\in \partial D_{\alpha_1} \cap  B(x_0,\, c_3\cdot \Delta)
\qquad\ \ \text{or}\qquad\ \ 
x^\beq
\in \partial D_{\alpha_2} \cap  B(x_0,\, c_3\cdot \Delta)\,,
\tag 3.37
$$
for some
$$
x_0\in G_{\!_{\alpha_1,\alpha_2}}\cap 
\bigcup_{\gamma \in \Cal T_{\alpha_2}}\text{Im}(\gamma)\,.
\tag 3.38
$$
Then, by (3.3), (3.24), (3.37), one has 
$$
\align
\big|x^\bem
-x_0\big|
&< (1 + c_3)\cdot \Delta < \delta_2\,,
\tag 3.39
\\
\noalign{\medskip}
\big|x^\beq
-x_0\big|
&< (1 + c_3)\cdot \Delta < \min\{\delta_1,\, \overline\lambda\,\}\,.
\tag 3.40
\endalign
$$
Thus, observing that by property $i)$ one has
$$
y^\bem(t)=x^g(t-\tau_1;\, 0,\, x^\bem)
\qquad\  \forall~t\in [\,\tau_0,\,\tau_1]\,,
\tag 3.41
$$
and because of $(3.25)_1$, (3.39),
applying property {\bf P2}
we deduce that 
there is another  backward trajectory
$t \mapsto x^g(t;\, 0,\, x_0),$ \ $t\in [-\widehat\tau_{x_0,x^\bem}^g,\,0],$ 
\, of (1.11)
arriving in \, $x_0,$ \, such that
%
$$
\gather
\big|\widehat\tau_{x_0,x^\bem}^g-\tau_1+\tau_0\big|
<c_2(1 + c_3)\cdot \Delta,
\tag 3.42
\\
\noalign{\medskip}
\big|y^\bem(t)- x^g\big(t-\tau_0-\widehat\tau_{x_0,x^\bem}^g;\, 0,\, x_0\big)\big|
< c_2 (1 + c_3)\cdot \Delta
\qquad\quad\forall~t\in~[\,\tau_0,\ \min\{\tau_0+\widehat\tau_{x_0,x^\bem}^g,\, \tau_1\}].
\tag 3.43
\endgather
$$
To fix the ideas assume that
$$
\tau_1-\tau_0<\widehat\tau_{x_0,x^\bem}^g<\tau_2-\tau_0\,.
\tag 3.44
$$
Then, 
observing that by property $ii)$ there holds (3.32),
using (3.9), (3.42)-(3.43), we obtain
$$
\align
\big|y^\beq(t)- x^g\big(t-\tau_0-\widehat\tau_{x_0,x^\bem}^g;\, 0,\, x_0\big)\big|
&\leq \big|y^\beq(t)-x^\beq\big|+\big|x^\beq-x^\bem\big|+\big|x^\bem-
x^g(\tau_1-\tau_0-\widehat\tau_{x_0,x^\bem}^g;\, 0,\, x_0)\big|+
\\
\noalign{\smallskip}
&\qquad +\big|x^g(\tau_1-\tau_0-\widehat\tau_{x_0,x^\bem}^g;\, 0,\, x_0)
- x^g(t-\tau_0-\widehat\tau_{x_0,x^\bem}^g;\, 0,\, x_0)\big|
\\
\noalign{\medskip}
&\leq 2M\cdot 
|t-\tau_1|+ \Delta
+  c_2 (1 + c_3)\cdot \Delta
\\
\noalign{\medskip}
&< 2(M+1)c_2 (1 + c_3)\cdot \Delta
\\
\noalign{\medskip}
&\hskip 0.9in \forall~t\in~]\tau_1,\, \tau_0+\widehat\tau_{x_0,x^\bem}^g]\,.
\tag 3.45
\endalign
$$
On the other hand, if \, $x^\beq\neq x_0,$ \, 
since 
(3.40) together with (3.20), imply (3.30),
using~$(3.25)_2$ 
we deduce as in Case~1 that (3.31) holds.
Hence, thanks to (3.40), 
by property {\bf P1}-{\it b)} it follows that
there exists \, $\widehat\tau_{x_0,x^\beq}^{\alpha_2}\in[0,\, t_{\alpha_2}^+(x_0)]$ 
\, such that
%
$$
\gather
x^{\alpha_2}\big(\widehat\tau_{x_0,x^\beq}^{\alpha_2};\, 0,\, x_0\big)
\in \partial D_{\alpha_2},
\tag 3.46
%
\\
\noalign{\medskip}
\big|\tau_2-\tau_1
-\widehat\tau_{x_0,x^\beq}^{\alpha_2}\big|
< c_1 (1 + c_3)\cdot \Delta\,.
\tag 3.47
\endgather
$$
Then, 
setting
$$
\sigma_{_{\bem,\beq}}\doteq
\tau_0+\widehat\tau_{x_0,x^\bem}^g+\widehat\tau_{x_0,x^\beq}^{\alpha_2}\,,
\tag 3.48
$$
and relying on (3.42), (3.47), we derive
$$
\big|\sigma_{_{\bem,\beq}}-\tau_2\big|
\leq \big|\tau_0+\widehat\tau_{x_0,x^\bem}^g-\tau_1\big|+
\big|\tau_1-\tau_2+\widehat\tau_{x_0,x^\beq}^{\alpha_2}\big|
<(c_1+c_2)(1 + c_3)\cdot \Delta\,,
\tag 3.49
$$
while, using (3.11), (3.40), (3.42), and because of (3.32),  we get
$$
\align
\big|y^\beq(t)-x^{\alpha_2}\big(t;\, \tau_0+\widehat\tau_{x_0,x^\bem}^g,\,x_0\big) \big|
&\leq c_0\cdot
\big(\big|x^\beq-x_0\big|+\big|\tau_1-\tau_0-\widehat\tau_{x_0,x^\bem}^g\big|\big)
\\
\noalign{\medskip}
&\leq c_0\cdot(1+c_2)(1 + c_3)\cdot \Delta
\\
\noalign{\medskip}
&\hskip 0.8in \forall~t\in~]\tau_0+\widehat\tau_{x_0,x^\bem}^g,\  
\min\{\tau_2,\, \sigma_{_{\bem,\beq}}\}]\,.
\tag 3.50
\endalign
$$
Thus, define
$$
\Phi_{_{\bem,\beq}}(t)\doteq
\cases
x^g\big(t-\tau_0-\widehat\tau_{x_0,x^\bem}^g;\, 0,\, x_0\big) \quad 
&\forall~ t\in ~[\,\tau_0,\, \tau_0+\widehat\tau_{x_0,x^\bem}^g]\,,
\\
\noalign{\medskip}
x^{\alpha_2}\big(t;\, \tau_0+\widehat\tau_{x_0,x^\bem}^g,\,x_0\big)
\quad &\forall~
t\in ~]\tau_0+\widehat\tau_{x_0,x^\bem}^g,\, \sigma_{_{\bem,\beq}}]\,,
\endcases
\tag 3.51
$$
and observe that, by construction and because of (3.46),
\, $\Phi_{_{\bem,\beq}}(\cdot)$ \, is a solution of (1.11) verifying~(3.6).
Moreover, from (3.43), (3.45), (3.50), and (3.47), 
it follows that  \, $\Phi_{_{\bem,\beq}}(\cdot)$ \,
satisfies 
the estimates (3.4), (3.5), (3.7) with the constant  
\, $C^{\text{\i v}}>2(M+c_0)(1 +c_2) (1 + c_3),$ \,
which shows that (3.51) provides the map
 required by the Proposition
whenever (3.37) holds.

\vs

\noindent
{\smbf 5.} \ 
{\smc Case~2.} \ Assume that 
$$
x^\bem
\in \partial D_{\alpha_1} \cap  B(x_0,\, c_3\cdot \Delta)
\qquad\ \ \text{for \ \ some}\qquad\ \ 
x_0\in G_{\!_{\alpha_1,\alpha_2}}\setminus
\bigcup_{\gamma \in \Cal T_{\alpha_2}}\text{Im}(\gamma)\,,
\tag 3.52
$$
and that
$$
x^\beq
\in \partial D_{\alpha_2} \cap  B(y_0,\, c_3\cdot \Delta)
\qquad\ \ \text{for \ \ some}\qquad\ \ 
y_0\in G_{\!_{\alpha_1,\alpha_2}}\setminus
\bigcup_{\gamma \in \Cal T_{\alpha_2}}\text{Im}(\gamma)\,.
\tag 3.53
$$
Observe that, by (3.3), (3.24), (3.52)-(3.53), one has
$$
\align
\big|x^\bem-x_0\big|
&< c_3\cdot \Delta < \delta_2\,,
\tag 3.54
\\
\noalign{\medskip}
\big|x^\beq-y_0\big|
&< c_3\cdot \Delta < \min\{\delta_1,\, \overline\lambda\,\}\,,
\tag 3.55
\\
\noalign{\medskip}
\big|x^\beq-x_0\big|
&< (1 + c_3)\cdot \Delta < \overline \lambda\,.
\tag 3.56
\endalign
$$
But then from (3.55)-(3.56), because of (3.20), we deduce that \, $x_0=y_0.$ \,
Hence, by (3.20), (3.52), (3.55), it follows that \, $x_0,\, x^\beq$ \,
belong to the same connected component of \, $\partial D_{\alpha_2}\setminus
\bigcup_{\gamma \in \Cal T_{\alpha_2}}
\text{Im}(\gamma).$ \, Moreover, since 
\, $y^\beq(\cdot)$ \, satisfies property $ii),$ we deduce also that (3.30)-(3.32)
hold. Therefore, relying on (3.55), and applying property {\bf P1}-{\it a)},
we derive
$$
\big|\tau_2-\tau_1-t_{\alpha_2}^+(x_0)\big|
< c_1 \,c_3\cdot \Delta\,.
\tag 3.57
$$
On the other hand, since by property $i)$ there holds (3.41), 
and because of $(3.25)_1$, (3.54),
applying property {\bf P2} as in Case~2
we deduce that 
there is another  backward trajectory
$t \mapsto x^g(t;\, 0,\, x_0),$ \ $t\in 
[-\widehat\tau_{x_0,x^\bem}^g,\,0],$ \, of (1.11)
arriving in \, $x_0,$ \, for which the estimates (3.42), (3.43), (3.45)
are verified.
Thus, setting
$$
\sigma_{_{\bem,\beq}}\doteq
\tau_0+\widehat\tau_{x_0,x^\bem}^g+t_{\alpha_2}^+(x_0)\,,
\tag 3.58
$$
and relying on (3.42), (3.57), we derive
$$
\big|\sigma_{_{\bem,\beq}}-\tau_2\big|
\leq \big|\tau_0+\widehat\tau_{x_0,x^\bem}^g-\tau_1\big|+
\big|\tau_1-\tau_2+\widehat\tau_{x_0,x^\beq}^{\alpha_2}\big|
<(c_1+c_2)(1 + c_3)\cdot \Delta\,,
\tag 3.59
$$
while, using (3.11), (3.42), (3.55), and because of (3.32),  we get
$$
\align
\big|y^\beq(t)-x^{\alpha_2}\big(t;\, \tau_0+\widehat\tau_{x_0,x^\bem}^g,\,x_0\big) \big|
&\leq c_0\cdot
\big(\big|x^\beq-x_0\big|+\big|\tau_1-\tau_0-\widehat\tau_{x_0,x^\bem}^g\big|\big)
\\
\noalign{\medskip}
&\leq c_0\cdot(1+c_2)(1 + c_3)\cdot \Delta
\\
\noalign{\medskip}
&\hskip 0.8in \forall~t\in~]\tau_0+\widehat\tau_{x_0,x^\bem}^g,\  
\min\{\tau_2,\, \sigma_{_{\bem,\beq}}\}]\,.
\tag 3.60
\endalign
$$
Observe now that the map \, $\Phi_{_{\bem,\beq}}(\cdot)$ \, 
defined in (3.51) is a solution of (1.11) verifying~(3.6)
since, by the definition (3.12) of the quantity \, $t_{\alpha_2}^+,$ \,
one has
$$
\Phi_{_{\bem,\beq}}(\sigma_{_{\bem,\beq}})=
x^{\alpha_2}\big(t_{\alpha_2}^+(x_0);\, 0,\, x_0\big)\in \partial D_{\alpha_2}\,.
$$
Moreover, from (3.43), (3.45), (3.60), and (3.57), 
it follows that  \, $\Phi_{_{\bem,\beq}}(\cdot)$ \,
satisfies 
the estimates (3.4), (3.5), (3.7) with the constant  
\, $C^{\text{\i v}}>2(M+c_0)(1 +c_2) (1 + c_3),$ \,
thus showing that (3.51) provides the map
 required by the Proposition
even in the case where (3.52)-(3.53) hold.
This completes the proof of Proposition~3.1 since, because of (3.1)-(3.2)
one has \, $x^\bem \in  \partial D_{\alpha_1}, \, x^\beq \in  
\partial D_{\alpha_2},$ \,
and hence the above three Cases~1-3 cover all the possibilities.
\fine
\vs

\v
\v
\n{\smbf Completion of the Proof of Theorem~1.} 
\v

Let \, $g$ \, be a 
uniformly
non-zero, polygonal  patchy
vector field on $\R^2,$ satisfying condition~{\bf (C)}.
Fix \, $T>0,$ and a
compact set $K\subset \Bbb R^2.$ 
Observe that, thanks to  Proposition~2.4, 
there exist contants \,
$C'''=C'''(K,T), 
\,  \delta'''=\delta'''(K,T)>0,$ 
so that,
for every BV perturbation \, $w=w(t)$ \,
with \, $\tv\{w\}< \delta''',$ \, 
and for every solution \, $y: [0,T]\mapsto\Bbb R^2$ \, of~(1.1),
with \, $y(0)\in K,$ \,  there is  a CCS of (1.11)
\, $y^\ds: [0,T]\mapsto\Bbb R^2,$ \, that satisfies  (2.17)-(2.18).
For this reason, in order to 
estabilish Theorem~1, it will be sufficient to 
take in consideration  only 
perturbed solutions  of (1.1) that are CCS of (1.11),
and derive for any such solution \, $y(\cdot)$ \, a
linear estimate of the distance from
some solution \, $x(\cdot)$ \, of (1.11), of the type
$$
\big\|x - y\big\|_{\L^\infty}\leq
C\cdot \Delta(y),
\tag 3.61
$$
with
$$
\Delta(y)\doteq 
\sum_{k=2}^{m}
\big|y(\tau_k^{\ +})-y(\tau_k)\big|.
\tag 3.62
$$
Here \, $0=\tau_1<\tau_2< \cdots < \tau_{m+1}=T$ \,
denotes the partition of \, $[0,\,T]$ \, defined as in (2.4)-(2.5),
that is
induced by \, $y(\cdot)$ \, according with the properties  $a''')$-$b''')$
of a CCS stated 
in Proposition~2.4. To this end we will establish the following 
\v

\n{\bf Lemma~3.2.} {\it 
%
In the same setting of Proposition~3.1,
given $T>0$ and any compact set $K \subset \Bbb R^2,$
there exist constants \,
$C^{\text{\sl v}}=C^{\text{\sl v}}(K,T), \, \delta^{\text{\sl v}}=
\delta^{\text{\sl v}}(K,T)>0,$ \,
for which the following hold.

For every CCS of (1.11) \, $y : [0,T]\mapsto\Bbb R^2,$ \, $y(0)\in K,$ \,
such that \, $\Delta(y)<\delta^{\text{\sl v}},$ \, 
letting \, $0=\tau_1<\tau_2< \cdots < \tau_{m+1}=T$ \,
be the partition of \, $[0,\,T]$ \, induced by \, $y(\cdot)$ \, 
according with (2.4)-(2.5),
there exist a sequence of CCS \, $y^k: [0,T]\mapsto \Bbb R^2,$ \,
$1\leq k < m',$ \, together with a 
sequence of points \,
$\widehat \tau_2,\,\widehat \tau_3, \cdots ,
\widehat \tau_{m'}=T\in ~]0,\,T],$ \,
having the following properties.
\v
\parindent 28pt

\item{\rm I)} \, $y^k\!\restriction_{[0,\,\widehat \tau_{k+1}]}$ \, is a 
Carath\'eodory solution of (1.11)
and, letting
$$
\alpha_k\doteq\max\big\{
\alpha~;~ \alpha\in \text{Im}\big((\alpha^* \circ y^k)\!
\restriction_{[0,\,\widehat \tau_{k+1}]}\!\big)
\big\}\,,
\tag 3.63
$$
one has
%
$$
\gathered
\alpha_k> \alpha_{k-1}
\qquad\text{if} \qquad k>1,
\\
\noalign{\medskip}
y^k(\,\widehat \tau_{\!_{k+1}})\in \partial D_{\alpha_k}
\qquad\text{if} \qquad k<m'-1,
\endgathered
\tag 3.64
$$
(where \, $D_{\alpha_k}$ \, denotes a polygonal
domain defined as in (2.3));
\v
\item{\rm II)} \, If \, $k<m'-1,$ \, one has
$$
\gather
y^k(\,\widehat \tau_{\!_{k+1}}^{\ \ \, +})\in \bigcup_{\gamma>\alpha_k} D_\gamma\,,
\tag 3.65
\\
\noalign{\medskip}
\alpha^*\big(y^k(t)\big)>\alpha_k \qquad\quad\forall~t\in 
~]\widehat \tau_{k+1},\,T]\,;
\tag 3.66
\endgather
$$

\v
\item{\rm III)} \, If \, $k>1,$ \, there holds
$$
\big|y^k(t)-y^{k-1}(t)\big| \leq C^{\text{\sl v}}\cdot \Delta(y^{k-1})
\qquad\quad\forall~t\in[0,\,T]\,,
\eqno (3.67)_k
$$
$$
\Delta(y^k)\leq\big(1+2 C^{\text{\sl v}}\,\big)^{k-1}\cdot \Delta(y)\,.
\eqno (3.68)_k
$$
}
\parindent 20pt
\vs

\noindent
{\bf Proof of Lemma~3.2.}  \ 
Fix \, $T>0,$ and a
compact set $K\subset \Bbb R^2.$ 
We may assume that
all CCS of~(1.11) \, $y : [0,T]\mapsto\Bbb R^2$ \, starting in a neighborhood
\, $B(K,\, \rho),$ \, $\rho>0,$ \, of  \, $K,$ \,
and having some uniform
bound~$\delta_4>0$ \, on the total amount of jumps~$\Delta(y)$, \,
take values in some compact 
set \, $K'.$ \,
Consider the set of indices
$$
\Cal A_{K'}\doteq \big\{\alpha \in \Cal A~~;~\Omega_\alpha\cap K'\neq \emptyset
\big\}\,,
$$
denote \, $N=|\Cal A_{K'}|$ \, the number of elements in \, $\Cal A_{K'},$ \,
and let \, $C^{\text{\sl \i v}}=C^{\text{\sl \i v}}(K',T), 
\, \delta^{\text{\sl \i v}}=
\delta^{\text{\sl \i v}}(K',T)>0,$ \, be constants chosen according 
with Proposition~3.1.
By the properties of a CCS and because of the regularity 
of the vector fields \, $g_\alpha$, \, 
$\alpha \in \Cal A_{K'}\,,$ \,
there will be some constant \, $c_4>0$ \, so that, for any CCS of (1.11) 
\, $y : [0,2T]\mapsto K',$ \, 
one has
$$
\big|y(t)-y(t')\big|\leq
c_4\cdot \big(|t-t'|+\Delta(y)\big) 
\qquad\quad \forall ~t,t'\in [0,\,2T]\,.
\tag 3.69
$$
Then, set
$$
\delta^{\text{\sl v}}\doteq
\min
\left\{
\frac{\delta^{\text{\sl \i v}}}{\big(1+2 C^{\text{\sl v}}\,\big)^N}\,,\
\frac{\delta_4}{\big(1+2 C^{\text{\sl v}}\,\big)^N}\,,\
\frac{\rho}{N \big(1+2 C^{\text{\sl v}}\,\big)^N}\,,\
\frac{T}{C^{\text{\sl \i v}}} 
\right\},
\qquad\quad C^{\text{\sl v}}\doteq 
3 c_4(1+2  C^{\text{\sl \i v}})\,, 
\tag 3.70
$$
and consider a CCS of 
(1.11) \, $y : [0,T]\mapsto\Bbb R^2,$ \, $y(0)\in K,$ \,
with
$$
\Delta(y) < \delta^{\text{\sl v}}.
\tag 3.71
$$
Let \, $0=\tau_1<\tau_2< \cdots < \tau_{m+1}=T$ \,
be the partition of \, $[0,\,T]$ \, induced by \, $y(\cdot)$ \, 
according with~(2.4)-(2.5).
We shall construct 
the 
sequence of of  CCS \, $y^k: [0,T]\mapsto \Bbb R^2$ \,
and of points \, $\widehat\tau_{k+1}\,,$ \, 
enjoing the properties I-III, 
applying  Proposition~3.1 and
proceeding by induction on \, $k\geq 1.$ \,
Set
$$
\widehat\tau_2 \doteq \tau_2, \qquad\qquad
y^1(y)\doteq y(t)\qquad\forall~t\in[0,T]\,.
\tag 3.72
$$
and, if \, $m=1,$ i.e. \, $\tau_2=T,$ \, set \, $m'\doteq 2,$ \, otherwise let \, $m'>2.$ \,
Observe that, by the properties  $a''')$-$b''')$
of a CCS stated 
in Proposition~2.4, the point \, $\widehat \tau_2$ \,
and the map \, $y^1(\cdot)$ \, in (3.72) clearly verify the conditions I-II of
Lemma~3.2.
Next, assume to have constructed, for some \, $1<k \leq N,$ \,   
a sequence of CCS \, $y^1, \dots , y^{k-1},$ \, 
together with a sequence of
points \, 
$0<\widehat\tau_2, \cdots ,\widehat\tau_k<T,$ \,
enjoing the 
properties~I-III, with \, $C^{\text{\sl v}}\,,$ \, as in (3.70).
Set
$$
\widehat \tau_{\!_{k+1}}'
\doteq
\sup \big\{t\in~]\,\widehat\tau_k,\,T]~;~ \alpha^*\big(y^{k-1}(t)\big)=
\alpha^*\big(y^{k-1}(\widehat\tau_k^{\  +})\big)
\big\}\,,
\tag 3.73
$$
and observe that, because of I-III,
and by $(3.68)_{k-1}$, (3.70)-(3.71), 
the maps
$$
y^{k,\bem}\doteq y^{k-1}\!\restriction_{[0,\,\widehat \tau_k]},
\qquad\quad 
y^{k,\beq}\doteq y^{k-1}\!
\restriction_{]\,\widehat \tau_k,\,\widehat \tau_{\!_{k+1}}']},
$$
have the properties $i)-iii)$ stated in Proposition~3.1.
Moreover, since $(3.67)_h$, $(3.68)_h$, $h=2, \dots , k-1,$
together with (3.70)-(3.71), imply
$$
\align
\big|y^{k-1}(0)-y(0)\big|
&\leq \sum_{h=2}^{k-1} \big|y^h(0)-y^{h-1}(0)\big|
\\
\noalign{\smallskip}
&\leq \sum_{h=2}^{k-1} \big(1+2 C^{\text{\sl v}}\,\big)^{h-2}\cdot \Delta(y)
\\
\noalign{\smallskip}
&\leq N \big(1+2 C^{\text{\sl v}}\,\big)^N \cdot \Delta(y)<\rho\,,
\tag 3.74
\endalign
$$
by the above assumptions it follows that \, $y^{k-1}(\cdot),$ \,
and hence \, $y^{k,\bem}(\cdot), \, y^{k,\beq}(\cdot),$ \,
take values in the set~$K'.$ \,
Then, letting \, 
$$
\Phi_{_{\bem,\beq}}^k\doteq\Phi\big[y^{k,\bem},\, y^{k,\beq}\big]:~
[0,\, \sigma_{_{\bem,\beq}}^k] \longmapsto \Bbb R^2
$$ \,
be the Carath\'eodory solution of (1.11)
provided by Proposition~3.1, and denoting by \, $x^{g}\big(t;\, t_0,\,x_0)$ \,
a Carath\'eodory solution of (1.11) starting 
from \, $x_0$ \, at time \, $t_0,$ \,
set
$$
\gather
\widehat \tau_{\!_{k+1}}\doteq \min\big\{\sigma_{_{\bem,\beq}}^k,\ T\big\}\,,
\qquad\quad
\widehat \tau_{\!_{k+1}}''\doteq 
\min\big\{T,\ T+\widehat \tau_{\!_{k+1}}-\widehat \tau_{\!_{k+1}}'\big\}\,,
\tag 3.75
\\
\noalign{\medskip}
y^k(t)\doteq
\cases
\Phi_{_{\bem,\beq}}^k(t) \quad\
&\text{if}\qquad t\in [0, \, \widehat \tau_{\!_{k+1}}]\,,
\\
\noalign{\smallskip}
y^{k-1}\big(t+\widehat \tau_{\!_{k+1}}'-\widehat \tau_{\!_{k+1}}\big)
\quad\
&\text{if}\qquad t\in ~]\,\widehat \tau_{\!_{k+1}},\, 
\widehat \tau_{\!_{k+1}}'']\,,
\\
\noalign{\smallskip}
x^g\big(t;\ \widehat \tau_{\!_{k+1}}'',\, y^{k-1}(T)\big)
\quad\
&\text{if}\qquad t\in ~]\,\widehat \tau_{\!_{k+1}}'',\, T]\,.
\endcases
\tag 3.76
\endgather
$$
Next, if \,  $\widehat \tau_{\!_{k+1}}=T,$ \,
set \, $m'\doteq k+1,$ \, otherwise let \,  \, $m'>k+1.$ \,
By construction, and because \, $\Phi_{_{\bem,\beq}}^k$ \,
satisfies condition (3.6) of Proposition~3.1, the map in (3.75)
defines a CCS of (1.11) that enjoys the properties~I-II.
Moreover,
by (3.4), (3.5), (3.7), and because of the above definition (3.75) of 
\, $\widehat \tau_{\!_{k+1}}$, \, we derive
$$
\gather
\big|\widehat \tau_{\!_{k+1}}-\widehat \tau_{\!_{k+1}}'\big|
\leq C^{\text{\sl v}}\cdot \Delta(y^{k-1})\,,
\tag 3.77
\\
\noalign{\medskip}
\big|y^k(t)-y^{k-1}(t)\big| \leq C^{\text{\sl v}}\cdot \Delta(y^{k-1})
\qquad\quad \forall~t\in
[0,\ \min\{\widehat \tau_{\!_{k+1}},\, \widehat \tau_{\!_{k+1}}'\}]\,.
\tag 3.78
\endgather
$$
On the other hand, using (3.69), (3.77)-(3.78), 
if \, $\widehat \tau'_{\!_{k+1}}<\widehat \tau_{\!_{k+1}}$ \, we obtain
$$
\align
\big|y^k(t)-y^{k-1}(t)\big| 
&\leq \big|\Phi_{_{\bem,\beq}}^k(t)-\Phi_{_{\bem,\beq}}^k(\widehat \tau_{\!_{k+1}}')\big|+
\big|y^k(\widehat \tau_{\!_{k+1}}')-
y^{k-1}(\widehat \tau_{\!_{k+1}}')\big|+
\big|y^{k-1}(t)-y^{k-1}(\widehat \tau_{\!_{k+1}}')\big|
\\
\noalign{\smallskip}
&\leq 2c_4 \cdot \big|t-\widehat \tau_{\!_{k+1}}'\big|+
c_4\cdot \big(\Delta(\Phi_{_{\bem,\beq}}^k)+ \Delta(y^{k-1})\big)+
C^{\text{\sl v}}\cdot \Delta(y^{k-1})
\\
\noalign{\smallskip}
&\leq 3 c_4(1+2  C^{\text{\sl v}})\cdot \Delta(y^{k-1})
\qquad\qquad \forall~t\in~ ]\,\widehat \tau_{\!_{k+1}}',\, \widehat \tau_{\!_{k+1}}]\,,
\tag 3.79
\endalign
$$
while, in the case \, $\widehat \tau_{\!_{k+1}}<T,$ \, we get
$$
\align
\big|y^k(t)-y^{k-1}(t)\big| 
&\leq c_4 \cdot\big(\big|\widehat \tau_{\!_{k+1}}-
\widehat \tau_{\!_{k+1}}'\big|+\Delta(y^{k-1})\big)
\\
\noalign{\smallskip}
&\leq c_4(1+C^{\text{\sl v}})\cdot \Delta(y^{k-1})
\qquad\qquad \forall~t\in~ ]\,\widehat \tau_{\!_{k+1}},\, T]\,.
\tag 3.80
\endalign
$$
From (3.78)-(3.80) we recover the estimate $(3.67)_k\,,$ 
with \, $C^{\text{\sl v}}$ \, as in (3.70), while $(3.68)_{k-1}\,,$
together with $(3.67)_k$, immediately yields
$(3.68)_k\,,$ showing that the map in (3.76) enjoys also the property~III.
\v

To complete the proof of Lemma~3.2, observe that proceeding by induction on
\, $k\geq 2,$ \, either we find some \, $m' \leq N$ \, 
such that \, $\widehat\tau_{m'}=T,$ \, or else we construct a 
sequence of CCS of (1.11)
\, $y^1, \dots , y^{N-1},$ \, together with a sequence of
points \, 
$0<\widehat\tau_2, \cdots ,\widehat\tau_N<T,$ \,
enjoing the properties~I-III. But then, if we
define \, $y^N$ \, and \, $\widehat \tau_{\!_{N+1}}$ \,
according with (3.75)-(3.76), we certainly find \, $\widehat\tau_{\!_{N+1}}=T,$ \,
since otherwise, relying on properties I-III, and performing
a computation as in (3.74), one deduces
$$
\gather
y^N(t)\in \Bbb R^2\setminus \bigcup_{\alpha\in \Cal A_{K'}} \Omega_\alpha
\qquad\ \
\forall~t\in~]\widehat\tau_{\!_{N+1}},\, T]\,,
\\
\noalign{\medskip}
\big|y^N(0)-y(0)\big|<\rho,
\endgather
$$
which contradicts the assumption that every 
CCS \, $y(\cdot)$ \, of~(1.11), starting  inside \, $B(K,\, \rho),$ \, 
and satisfying the bound \, $\Delta(y)<\delta_4$,  \,
takes values in the set \, $\bigcup_{\alpha\in \Cal A_{K'}} \Omega_\alpha.$\,
This concludes the proof of the Lemma.
\fine
\vs

We are in the position now to
complete the proof of Theorem~1,
relying on Lemma~3.2. Let \, 
\, $C^{\text{\sl v}}=C^{\text{\sl v}}(K,T), \ \delta^{\text{\sl v}}=
\delta^{\text{\sl v}}(K,T)>0$ \, be constants chosen according with Lemma~3.2
and consider a  CCS of (1.11) \, $y : [0,T]\mapsto\Bbb R^2,$ \, $y(0)\in K,$ \,
such that \, $\Delta(y)<\delta^{\text{\sl v}}.$ \, 
Then, letting \, $y^1=y, \, y^2, \, \dots y^{m'-1}$ \, 
be the 
sequence of CCS
provided by Lemma~3.2 that
enjoy the properties I-III, and using $(3.67)_k$-$(3.68)_k\,,$ we derive
$$
\align
\big|y^{m'-1}(0)-y(0)\big|
&\leq \sum_{k=2}^{m'-1} \big|y^k(t)-y^{k-1}(t)\big|
\\
\noalign{\smallskip}
&\leq N \big(1+2 C^{\text{\sl v}}\,\big)^N\,.
\tag 3.81
\endalign
$$
By property I, \, $y^{m'-1} : [0, T] \to \Bbb R^2$ \, is a 
Carath\'eodory solution of (1.11), and hence (3.81)
yields the estimate in (3.61), taking \, 
$C=N \big(1+2 C^{\text{\sl v}}\,\big)^N,$ \,
which concludes the proof of the theorem.
\fine
\vsk

\n{\medbf 4 - Appendix.}
\vs

We provide here the proofs of the Lemmas~2.2-2.3
stated in Section~2. To this end we shall establish first
the following

\v
\n{\bf Lemma~4.1.} {\it Let $g$ be a 
uniformly
non-zero polygonal  patchy
vector field on $\R^2,$
associated to a family of polygonal patches
$\big\{ (\Omega_\alpha,~g_\alpha)  ; \alpha\in\Cal A\big\}$.
Assume that condition {\bf (C)} (stated in Section~2)
is satisfied.
Then, given $T,\, C>0$ and any compact set \, $K \subset \Bbb R^2,$ \,
there exist constants \, $C^{\text{\sl v\i}}=
C^{\text{\sl v\i}}(K,T,C)\geq C,$ \ 
$C^{\text{\sl v\i\i}}=C^{\text{\sl v\i\i }}(K,T,C), \
\delta^{\text{\sl v\i}}=\delta^{\text{\sl v\i}}(K,T,C)>0,$ \,
so that
the following property holds.

For every BV perturbation \, $w=w(t)$ \,
with \, $\tv\{w\}<\delta^{\text{\sl v\i}},$ \, 
and for every left continuous  solution 
\, $y: [0,T]\mapsto\Bbb R^2$ \, of~(1.1),
starting at some point \, $y_0\in K,$ \, 
for which the map \, $t\mapsto \alpha^*\big(y(t)\big)$ \, is
non-decreasing, 
letting
$$
\big\{\,
{\alpha_{i_1}},\dots,{\alpha_{i_m}}
\,\big\}
= 
\text{Im}\big(\alpha^* \circ y\big),
\tag 4.1
$$
with
$$
\alpha_{i_1}< \dots < \alpha_{i_m},
\tag 4.2
$$
and setting
$$
\gather
D_{\alpha_{i_j}}
\doteq
\Omega_{\alpha_{i_j}} \setminus 
\displaystyle{\bigcup_{\beta > {\alpha_{i_j}}} \Omega_\beta}
\qquad\quad
j=1,\dots,m\,,
\tag 4.3
\\
\noalign{\smallskip}
]\tau_j,\,\tau_{j+1}]
\doteq
\Big\{
t\in [0,\,T]~:~y(t)\in 
 D_{\alpha_{i_j}}
\Big\}
\qquad\quad
j=1,\dots,m\,,
\endgather
$$
\v
\noindent
there holds
$$
\text{\rm meas}\,
\left(\bigcup_j
\Big\{
t\in ]\tau_j,\,\tau_{j+1}]~:~
d(y(t),\,\partial D_{\alpha_{i_j}})< 
C^{\text{\sl v\i}} \cdot \text{\rm Tot.Var.}\{w\}
\Big\}\right)
< C^{\text{\sl v\i\i }} \cdot \text{\rm Tot.Var.}\{w\}.
\tag 4.4
$$
}
\v

\noindent
{\bf Proof.} \ 

\n{\smbf 1.} \
Fix $T,\, C>0$ 
and a compact set \, $K\subset \Bbb R^2.$ \,
The smoothness
of the vector fields \, $g_\alpha$ \, guarantees that,
for BV perturbations \, $w=w(t)$ \, having some uniform bound
on the total variation,
every solution  \,
 $y : [0,T]\mapsto\Bbb R^2$ \, of (1.1),
starting at any point \, $y_0\in K,$ \, takes values in some compact
set \, $K'\subset \Bbb R^2.$ \,
Let
$$
\big\{
\Om_{\alpha_i}~:~i=1,\dots,N
\big\}
$$
be the finite collection of 
domains that intersect $K,$
with the indices ordered so that 
$$
\alpha_1< \dots < \alpha_{\!_N}.
\tag 4.5
$$
%
Let \, $D_{\alpha_i}, 
\ i=1,\dots,N,$ \, denote the
domains defined as in (4.3).
For each $\alpha,$ call \, 
$E_\alpha^{1}
, \dots ,E_\alpha^{^{p_\alpha}}
,$ \, 
and \, 
$r_\alpha^{1},
\dots ,r_\alpha^{^{p_\alpha}}
,$ \, 
respectively,
the edges of \, $D_{\alpha}$ \,
that form the boundary~$\partial D_\alpha,$ \, 
and the corresponding lines in which the edges are contained.
By construction, every edge \, $E_\alpha^{\ell}
$ \, 
is a part of
the boundary of  some~$\Om_\beta, \ \beta\geq\alpha.$ \,
Call \, $\bn_{{\boldsymbol\alpha}}^
{{\boldsymbol \ell}
}$ \,  
the  normal to~$E_\alpha^{\ell}
$ \, pointing towards the interior of \, $D_\alpha,$ \,
and
let $\varphi_\alpha^\ell(x)$ \,
denote the signed distance 
of the point $x$ from $r_\alpha^{\ell},$ i.e.
$$
\varphi_\alpha^\ell(x)\doteq
\cases
d\big(x,\ r_\alpha^{\ell}\big)
\quad &\text{if}\qquad x\in
r_\alpha^{\ell}+
\big\{\lambda\, \bn_{{\boldsymbol\alpha}}^
{{\boldsymbol \ell}
}
~:~ \lambda \geq 0\big\},
\\
\noalign{\smallskip}
-d\big(x,\ r_\alpha^{\ell}\big)
\quad &\text{if}\qquad x\in
r_\alpha^{\ell}+
\big\{\lambda\, \bn_{{\boldsymbol\alpha}}^
{{\boldsymbol \ell}
}
~:~ \lambda \leq 0\big\}.
\endcases
\tag 4.6
$$
Given any BV perturbation \, $w~:~[t_0,\,t_1]\mapsto D_{\alpha}$ \,
defined on some interval \, $[t_0,\,t_1]\subset [0,\,T],$ \,
and any solution \, $y~:~[t_0,\,t_1]\mapsto D_{\alpha}$ \,
of the Cauchy problem (1.1)-(1.2), with $y_0\in K,$
consider the map \, 
$\varphi_\alpha^\ell \circ y~:~[t_0,\,t_1]\mapsto \Bbb R.$ \,
One can easily verify that, 
for every Borel set \, $E\subset [t_0,\,t_1],$ \, the 
Radon measure \,
$\mu\doteq D \,(\varphi_{\alpha}^\ell\circ y)$ \,
satisfies 
$$
\bigg|
\mu(E)-
\int_E \big\langle \nabla \varphi_{\alpha}^\ell(y(t)),\ 
g_\alpha(y(t))
\big\rangle~dt 
\bigg|
\leq c_5\cdot
\tv\{w\},
\tag 4.7
$$
for some constant \, $c_5>0$ \, depending only on the 
compact set $K'$ and on the time interval $[0,T].$
Then, fix 
$$
C^{\text{\sl v\i}}>\max\{C,\,c_5\},
\tag 4.8
$$
and take
\, $\delta_5>0$ \, so that one has
$$
\Big\{
x\in D_{\alpha_i}
~:~ d(x,\,\partial D_{\alpha_i})> 2C^{\text{\sl v\i}} \cdot \delta_5
\Big\}
\neq \emptyset
\qquad\quad \forall~i.
\tag 4.9
$$
%
%
%
%
%
%
%
%
%
%
\v

Observe now that, since $g_\alpha$ are smooth, non-zero vector fields,
the transversality condition {\bf (C)} and the inward-pointing
condition (1.15) guarantee that there exists
some constant \, $c_6>0$ \, such that, for every \,
$\alpha_i, \ i=1, \dots ,N,$ \, and for every
\, $\ell=1,\dots,p_{\alpha_i},$ \, 
one of the  following two conditions
holds
$$
\alignat 4
\big \langle
g_{\alpha_i}(x),\ 
\bn_{{\boldsymbol\alpha}_{\bold i}}^{{\boldsymbol \ell}
}
\big\rangle
&\geq c_6
\qquad\quad &\forall~x \in 
D_{\alpha_i}\cap r_{\alpha_i}^{\ell}
\cap K',
\tag 4.10
\\
\noalign{\medskip}
\big \langle
g_{\alpha_i}(x),\ 
\bn_{{\boldsymbol\alpha}_{\bold i}}^{{\boldsymbol \ell}
}
\big\rangle
&\leq -c_6
\qquad\quad &\forall~x \in 
D_{\alpha_i}\cap r_{\alpha_i}^{\ell}
\cap K'.
\tag 4.11
\endalignat
$$
For each $\alpha_i,$ define the sets 
\, $\Cal I_{\alpha_i}$ \, and \, $\Cal O_{\alpha_i}$ \,
of \, (incoming and outgoing) \, indices
$$
\aligned
\Cal I_{\alpha_i}
&\doteq
\big\{1\leq \ell \leq p_{\alpha_i}~:~ (4.10) \quad\text{holds}
\big\},
\\
\noalign{\medskip}
\Cal O_{\alpha_i}
&\doteq
\big\{1\leq \ell \leq p_{\alpha_i}~:~ (4.11) \quad\text{holds}
\big\}.
\endaligned
\tag 4.12
$$
Because of the regularity assumptions on the fields $g_\alpha,$
there will be some constants
\, $0<\delta_6\leq \delta_5,$ \, $c_7>0,$ \, so that
$$
\aligned
&\sup
\Big\{
\big\langle g_{\alpha_i}(x),
\ \bn_{{\boldsymbol\alpha}_{\bold i}}^{{\boldsymbol \ell}
}
\big\rangle 
~:~\
x \in D_{\alpha_i}\cap
\Big(
r_{\alpha_i}^{\ell}
+ 
\big\{\lambda\, \bn_{{\boldsymbol\alpha}_{\bold i}}^
{{\boldsymbol \ell}
}
~:~ |\lambda| \leq 2 C^{\text{\sl v\i}}\cdot \delta_6\big\}
\Big)
\cap K',\qquad
\\
\noalign{\smallskip}
&\hskip 3in
i=1,\dots,N,\qquad
\ell \in \Cal I_{\alpha_i}
\Big\}
\geq c_7,
\\
\noalign{\bigskip}
&\sup
\Big\{
\big\langle g_{\alpha_i}(x),
\ \bn_{{\boldsymbol\alpha}_{\bold i}}^{{\boldsymbol \ell}
}
\big\rangle 
~:~\
x \in D_{\alpha_i}\cap
\Big(
r_{\alpha_i}^{\ell}
+ 
\big\{\lambda\, \bn_{{\boldsymbol\alpha}_{\bold i}}^
{{\boldsymbol \ell}
}
~:~ |\lambda| \leq 2 C^{\text{\sl v\i}}\cdot \delta_6\big\}
\Big)
\cap K',
\\
\noalign{\smallskip}
&\hskip 3in
i=1,\dots,N,\qquad
\ell \in \Cal O_{\alpha_i}
\Big\}
\leq -c_6.
\endaligned
\tag 4.13
$$
\vs

\noindent
{\smbf 2.} \ Consider now a
BV perturbation \, $w=w(t)$ \,
with \, $\tv\{w\}<\delta_6,$ \, 
and let \, $y: [0,T]\mapsto\Bbb R^2$ \, be a solution of~(1.1),
starting at some point \, $y_0\in K,$ \, for which there
is a partition \, $\tau_1<\tau_2<\cdots<\tau_{m+1}$ \, of \, $[0,\,T],$ \,
such that
$$
\alpha^*(y(t))=
\alpha_{i_j}
\qquad\quad \forall~t\in~]\tau_j, \tau_{j+1}],
\qquad j=1,\dots,m,
\tag 4.14
$$
with
$$
\alpha_{i_1}<\cdots<\alpha_{i_m}.
\tag 4.15
$$
For any \, $\alpha \in \{\alpha_{i_1},\dots,\alpha_{i_m}\},$ \
$\ell\in\{1,\dots,p_\alpha\},$ \,   call 
\, $S_{\alpha}^{\ell}$ \,
the connected component of the set
$$
\Big\{
x\in D_{\alpha}
~:~ d(x,\,\partial D_{\alpha})< C^{\text{\sl v\i}} \cdot \tv\{w\}
\Big\}
\bigcap
\Big(
r_\alpha^{\ell}
+ 
\big\{\lambda\, \bn_{{\boldsymbol\alpha}}^
{{\boldsymbol \ell}
}
~:~ 0\leq \lambda < C^{\text{\sl v\i}}\cdot \tv\{w\}\big\}
\Big)
\tag 4.16
$$
whose boundary contains the edge \, $E_\alpha^{\ell}
.$ \, 
Then, by similar computations as those used in \ABfs\
to establish
[2, \, Proposition~2.2],
relying on (4.7), (4.8), (4.13), 
one can establish the following claims.
\v

\noindent
{\bf Claim~1.} \ {\it If, for some \, 
$\alpha_{i_j}, \ j=1,\dots,m,$ \
$\ell\in \Cal I_{\alpha_{i_j}},$ \, 
and for some constant \, $0\leq c \leq 2 C^{\text{\sl v\i}},$ \,
there exists
 \, $t'\in ]\tau_j,\, \tau_{j+1}],$ \, such that
$$
\varphi_{\alpha_{i_j}}^\ell(y(t'\,))\geq 
c\cdot \tv\{w\},
\tag 4.17
$$
then, there holds
$$
\varphi_{\alpha_{i_j}}^\ell(y(t))> c-C^{\text{\sl v\i}}\cdot \tv\{w\}
\qquad\quad\forall~t\in [t',\,\tau_j].
\tag 4.18
$$
}
\v

\noindent
{\bf Claim~2.} \ {\it If, for some \, 
$\alpha_{i_j}, \ j=1,\dots,m,$ \
$\ell\in \Cal O_{\alpha_{i_j}},$ \, 
and for some constant \, $|c| \leq C^{\text{\sl v\i}},$ \,
there exists
 \, $t'\in ]\tau_j,\, \tau_{j+1}],$ \, such that
$$
\varphi_{\alpha_{i_j}}^\ell(y(t'\,))\leq 
c\cdot \tv\{w\},
\tag 4.19
$$
then, there holds
$$
\varphi_{\alpha_{i_j}}^\ell(y(t))< c+C^{\text{\sl v\i}}\cdot \tv\{w\}
\qquad\quad\forall~t\in [t',\,\tau_j].
\tag 4.20
$$
}
\v

\noindent
{\bf Claim~3.} \ {\it If, for some interval \,
$[t_1,\,t_2]\subset ]\tau_j,\, \tau_{j+1}],$ \
$j=1,\dots,m,$ \,
and for some $\ell\in \Cal I_{\alpha_{i_j}},$ \, 
the following two conditions hold
$$
\gather
\varphi_{\alpha_{i_j}}^\ell(y(t))
<2 C^{\text{\sl v\i}}\cdot \tv\{w\}
\qquad\forall~t\in[t_1,\,t_2],
\tag 4.21
\\
\noalign{\medskip}
\meas
\big\{t\in[t_1,\,t_2]~:~
y(t)\in S_{\alpha_{i_j}}^{\ell}
\big\}
> c_8\cdot \tv\{w\},
\tag 4.22
\endgather
$$
with
$$
c_8\doteq\frac{2 C^{\text{\sl v\i}}}{c_7},
\tag 4.23
$$
then, one has
$$
\varphi_{\alpha_{i_j}}^\ell(y(t_2))> C^{\text{\sl v\i}}\cdot \tv\{w\}.
\tag 4.24
$$
}
\v

\noindent
{\bf Claim~4.} \ {\it If, for some interval \,
$[t_1,\,t_2]\subset ~]\tau_j,\, \tau_{j+1}],$ \
$j=1,\dots,m,$ \,
and for some $\ell\in \Cal O_{\alpha_{i_j}},$ \, 
there holds
$$
\varphi_{\alpha_{i_j}}^\ell(y(t))
>-C^{\text{\sl v\i}}\cdot \tv\{w\}
\qquad\forall~t\in[t_1,\,t_2],
\tag 4.25
$$
together with the condition (4.22),
then, one has
$$
\varphi_{\alpha_{i_j}}^\ell(y(t_2))<0.
\tag 4.26
$$
}
\v

\noindent
From  Claim~1 (taking the
constant $c=2C^{\text{\sl v\i}}$), Claim~2 
(taking the constant $c=-C^{\text{\sl v\i}}$), and Claims~3-4,
it clearly follows that, for
every fixed \,
$j=1,\dots,m,$ \,
and every \, $\ell=1,\dots, p_{\alpha_{i_j}},$ \, 
there holds
$$
\meas
\big\{t\in]\tau_j,\, \tau_{j+1}]~:~
y(t)\in S_{\alpha_{i_j}}^{\ell}
\big\}
\leq c_8\cdot \tv\{w\}.
\tag 4.27
$$
Thus, observing that by construction we have
$$
\bigcup_j
\Big\{
t\in ]\tau_j,\,\tau_{j+1}]~:~
d(y(t),\,\partial D_{\alpha_{i_j}})< C^{\text{\sl v\i}} 
\cdot \text{\rm Tot.Var.}\{w\}
\Big\}=
\bigcup_{j,\,\ell}
S_{\alpha_{i_j}}^{\ell},
\tag 4.28
$$
from (4.27) we derive the estimate (4.4) with 
$$
C^{\text{\sl v\i\i }}\doteq \Big(\sum_{i=1}^N p_{\alpha_i}\Big)\cdot c_8,
\tag 4.29
$$
where \, $p_{\alpha_i}$ \, denotes the number of edges of 
the domain \, $D_{\alpha_i}$ \, defined as in (4.3), while $c_8$
is the constant defined in (4.23).
This completes the proof of Lemma~4.1, taking \, 
$\delta^{\text{\sl v\i}}=\delta_6$
and \, $C^{\text{\sl v\i}}$ \, as in~(4.8).
\fine
\vs

\noindent
{\bf Proof of Lemma~2.2.}  \ 
\v

\noindent 
{\smbf 1.} \ 
Fix $T>0$ and a
compact set $K\subset \Bbb R^2.$ 
Observe that, thanks to Proposition~2.1,
in order to estabilish Lemma~2.2
it will be sufficient to show that
there exist constants \, $C''=C''(K,T),$ \,  $\delta''=\delta''(K,T)>0,$ \,
so that, given any BV function \, $w(\cdot)$ \, with
\, $\tv\{w\}< \delta''$,\ ,
for every
left continuous solutions \, $y(\cdot)$ \,  of (1.1)
for which the map \, $t \mapsto \alpha^*\big(y(t)\big)$ \, 
is non-decreasing, there exists a
piecewise smooth function \,
$y^\di:[0,T]\mapsto\Bbb R^2$ \, enjoing the properties
$a')$-$c')$.

As in the proof of Lemma~4.1 we may assume that,
for BV perturbations \, $w=w(t)$ \, having some uniform
bound on the total variation,
every solution 
\, $y : [0,T]\mapsto\Bbb R^2$ \, of~(1.1), 
starting at a  point \, $y_0\in K,$ \, is contained in some
compact set \, $K'\subset \Bbb R^2.$ \ Let \,
$\big\{
\Om_{\alpha_i}~:~i=1,\dots,N
\big\}$ \, be  the collection of polygonal domains
that intersect  \, $K'.$ \,
We shall assume that every vector field \,$ g_{\alpha_i}$ \, is
defined on a  neighborhood 
\, $B(\Om_{\alpha_i},\, \rho), \ \rho>0,$ \,
of the domain \, $\Om_{\alpha_i}$ \, and,
for any fixed \, $t_0>0,$
\ $x_0\in B(\Om_{\alpha_i},\, \rho),$ \,
we will denote
by \, $t \mapsto x^{\alpha_i}(t;\, t_0,\, x_0)$ \
the 
solution of the Cauchy 
problem
%
$$
\dot y = g_{\alpha_i}(y),
\qquad
y(t_0)=x_0\,,
\tag 4.30
$$
and set
$$
M\doteq 
\sup\big\{|g_{\alpha_i}(y)|~:~y\in B(\Om_{\alpha_i},\, \rho),
\quad i=1,\dots,N\big\}\,.
\tag 4.31
$$
Similarly, for every given \, $w\in\BV,$ with \, $\tv\{w\}<\delta,$ \,
we denote by
\, $t \mapsto z^{\alpha_i}(t;\, w,\, t_0,\, x_0)$ \
the 
left-continuous solution of 
%
$$
\dot z = g_{\alpha_i}(z)+\dot w, 
\qquad
z(t_0)=x_0.
\tag 4.32
$$
For every \, $x_0\in B(\Om_{\alpha_i},\, \rho),$
\, $t_0>0,$ \,
we let \, $t^{\alpha_i,+}(t_0, x_0), \, t^{\alpha_i,-}(t_0, x_0)$ \,
denote
the time that is necessary to
reach the set \, 
$B(\Om_{\alpha_i},\, \rho)\setminus  D_{\alpha_i},$  \,
starting from $x_0$ at time
$t_0,$ and following, respectively,  the forward
and backward flow of $g_{\alpha_i},$ i.e. 
$$
\aligned
t_{\alpha_i}^{\,+}(t_0, x_0)
& \doteq
\inf\big\{t>t_0~:~ x^{\alpha_i}(t;\, t_0,\, x_0)
\in B(\Om_{\alpha_i},\, \rho)\setminus  D_{\alpha_i}
\big\}\,,
\\
\noalign{\smallskip}
t_{\alpha_i}^{\,-}(t_0, x_0)
& \doteq
\sup\big\{t<t_0~:~ x^{\alpha_i}(t;\, t_0,\, x_0)
\in B(\Om_{\alpha_i},\, \rho)\setminus  D_{\alpha_i}
\big\}\,.
\endaligned
\tag 4.33
$$
Since \,$ g_{\alpha_i}$ \, are smooth vector fields, and because of
the linear estimate (1.6), the Cauchy Problems 
(4.30), (4.32) are well posed.
Hence, there will be some constant $c_9>0$
so that there holds
$$
\qquad\quad
\Big|
x^{\alpha_i}(t;\, t_0,\, x_0)-
z^{\alpha_i}(t';\, w,\, t_1,\, x_1)
\Big|\leq c_9
\Big\{|t-t'|+|t_0-t_1|+|x_0-x_1|+\tv\{w\}\Big\}
\tag 4.34
$$
for any \, $t, \, t',\,  t_0,\, t_1,\, x_0, \, x_1, \,
\, w,$ \, and for every $\alpha_i.$ 
Moreover, recalling that \, $g_{\alpha_i}$ \, are uniformly non-zero 
vector fields
that satisfy the  inward-pointing
condition (1.15) and
the transversality condition {\bf (C)}, we deduce that
there will be  constants \, $c_{10}>1,\, \delta_7>0,$ \
so that, if 
$$
d\Big(x_0,\ B(\Om_{\alpha_i},\, \rho)\setminus  D_{\alpha_i}\Big)
<\delta_7\,,
\tag 4.35
$$
then
$$
\big|
t_{\alpha_i}^{\,\pm}(t_0, x_0)-t_0\big|
\leq c_{10} \cdot
d\Big(x_0,\ B(\Om_{\alpha_i},\, \rho)\setminus  D_{\alpha_i}\Big)\,.
\tag 4.36
$$
\v

\noindent
Set
$$
c_{11}\doteq 4(1+c_9)^2\cdot (1+ 2Mc_{10})^{N+1}\,,
\tag 4.37
$$
and let \, $C^{\text{\sl v\i}}=C^{\text{\sl v\i}}(K,T,2 c_{11})
\geq 2 c_{11},$ \ $C^{\text{\sl v\i\i}}=C^{\text{\sl v\i\i}}(K,T,2 c_{11}),$ \
$\delta^{\text{\sl v\i}}=\delta^{\text{\sl v\i}}(K,T,2 c_{11})>0,$ \, be 
constants chosen
according with Lemma~4.1. 
Then, fix any  \,
$w\in\BV$ with 
$$
\tv\{w\}<
\delta''\doteq \min\left\{\delta^{\text{\sl v\i}},\ \frac{\delta_7}{c_{11}}\,,\,
\frac{\rho}{c_{11}}
\right\},
\tag 4.38
$$
and consider a left continuous 
solution \, $y: [0,T]\mapsto\Bbb R^2$ \, of~(1.1),
starting at some point \, $y_0\in K,$ \linebreak  
for which the map \, $t\mapsto \alpha^*\big(y(t)\big)$ \, is
non-decreasing.
As an intermediate step towards the construction of the map \,
$y^\di:[0,T]\mapsto\Bbb R^2$ \, enjoing properties
$a')$-$c'),$ we shall first produce a piecewise smooth function
$\widetilde y:[0,T]\mapsto\Bbb R^2$ whose $\L^\infty$~distance
from \, $y(\cdot)$ \, is bounded by \, $c_{11}\cdot\tv\{w\},$ \,
and for which there is a partition \,
\, $0=\tau_1'<\tau_2'<\cdots<\tau_{m'+1}'=T$ \, 
of \, $[0, T],$ together with an increasing sequence of indices
\, $\alpha_{i'_1}< \cdots < \alpha_{i'_{m'}},$ \, 
so that
$\widetilde y(\cdot)$ is a classical solution of \, 
$\dot y = g_{\alpha_{i'_k}}(y)$ \,
on every interval \, 
$]\tau'_k, \, \tau'_{k+1}[\,,$ $1 \leq k \leq m',$ \,
but does not satisfy the condition 
\, $\widetilde y(t)\in D_{\alpha_{i'_k}}$ \,
for all \, $t\in ~]\tau'_k, \, \tau'_{k+1}[\,$.
\vs

\noindent
{\smbf 2.} \
In order to define the map \, $\widetilde y(\cdot),$ \, 
in connection with the partition 
\, $0=\tau_1<\tau_2<\cdots<\tau_{m+1}=T$ \, 
of \, $[0, T]$ \, induced by \, $y(\cdot)$ \,
according with (4.14)-(4.15),
consider the sequence of points \, $\tau_1'\doteq
0<\tau_2'<\cdots<\tau_{m'+1}'=T,$ \,
$m'\leq m,$ \, 
and of subindices \, 
$j(1)=1,\, j(2),\, \dots,\, j(m') \in \{1, \dots ,m\},$ \,
recursively defined as follows. 
If \, $\tau_2+ t_{\alpha_{i_1}}^{\,+}
\big(0,\, x^{\alpha_{i_1}}\big(\tau_2;\, 0,\, y(0^{\, +})\big)\big)\geq T,$
set  \, $\tau_2'\doteq T,$ \,
otherwise
set
$$
\tau_2'\doteq 
\tau_2+ t_{\alpha_{i_1}}^{\,+}
\big(0,\,  
x^{\alpha_{i_1}}\big(\tau_2;\, 0,\, y(0^{\, +})\big)\big)\,.
\tag 4.39
$$
Next, for all \, $1<k\leq m$ \, such that
\, $\tau'_k<T,$ \, 
let \, ${j(k)}$ \, be the subindex of \, $\alpha_i$ \,
for which there holds \, $\alpha_{i_{\!_{j(k)}}}=
\alpha^*\big(y({\tau'_k}^{\, +})\big),$ \,
so that
$$
\tau'_k\in [\tau_{_{j(k)}},\, \tau_{_{j(k)+1}}]\,,
\tag 4.40
$$
and set 
$$
\tau_{k+1}'\doteq
\min\Big\{T,\
\tau_{\!_{j(k)+1}}
\!-t_{\alpha_{i_{\!_{j(k)}}}}^{\,-}
\!\big(0,\,  y({\tau'_k}^{\, +})\big)+t_{\alpha_{i_{\!_{j(k)}}}}^{\,+}
\!\big(0,\,  
x^{\alpha_{i_{\!_{j(k)}}}}\big(\tau_{_{j(k)+1}};\, \tau_k',\, 
y({\tau'_k}^{\, +})\big)\big)\Big\}\,.
\tag 4.41
$$
Then, set 
$$
m'\doteq \max\big\{1\leq k\leq m~:~ \tau'_k<T\big\},
\qquad\quad
\tau'_{m'+1}\doteq T\,.
$$
Observe that, using (4.34), and because of (4.37)-(4.38), one finds,
for every \, $k\geq 1,$ 
$$
\align
d\Big(
x^{\alpha_{i_{\!_{j(k)}}}}&\big(\tau_{_{j(k)+1}};\, \tau_k',\, 
y({\tau'_k}^{\, +})\big),\ 
B(\Om_{\alpha_{i_{\!_{j(k)}}}},\, \rho)\setminus 
D_{\alpha_{i_{\!_{j(k)}}}}\Big)
\leq
\\
\noalign{\smallskip}
&\leq
\Big(
\big|
x^{\alpha_{i_{\!_{j(k)}}}}\big(\tau_{_{j(k)+1}};\, \tau_k',\, 
y({\tau'_k}^{\, +})\big)-
z^{\alpha_{i_{\!_{j(k)}}}}\big(\tau_{_{j(k)+1}};\, w,\,\tau_k',\, 
y({\tau'_k}^{\, +})\big)
\big|+
\\
\noalign{\smallskip}
&\qquad\quad  +
\big|
z^{\alpha_{i_{\!_{j(k)}}}}\big(\tau_{_{j(k)+1}};\, w,\,\tau_k',\, 
y({\tau'_k}^{\, +})\big)- y({\tau_{\!_{j(k)+1}}}^{\!+})
\big|
\Big)
\\
\noalign{\smallskip}
&\leq 
c_9\cdot \tv\{w\} +
\big|y\big(\tau_{_{j(k)+1}}\big)- y\big({\tau_{\!_{j(k)+1}}}^{\!+}\big)
\big|
\\
\noalign{\smallskip}
&\leq (1+c_9)\cdot \tv\{w\}< \delta_7\,.
\tag 4.42
\endalign
$$
Thanks to (4.42), and because of (4.35), we can apply (4.36) obtaining
$$
\align
\big|
t_{\alpha_{i_{\!_{j(k)}}}}^{\,+}
&\big(0,\,  
x^{\alpha_{i_{\!_{j(k)}}}}\big(\tau_{_{j(k)+1}};\, \tau_k',\, 
y({\tau'_k}^{\, +})\big)\big)
\big|\leq 
\\
\noalign{\smallskip}
&\leq
c_{10}\cdot d\Big(\!
x^{\alpha_{i_{\!_{j(k)}}}}\!\big(\tau_{_{j(k)+1}};\, \tau_k',\, 
y({\tau'_k}^{\, +})\big),\ 
B(\Om_{\alpha_{i_{\!_{j(k)}}}},\, \rho)\setminus 
D_{\alpha_{i_{\!_{j(k)}}}}\Big)
\\
\noalign{\smallskip}
&\leq c_{10}\cdot\Big(
c_9\cdot \tv\{w\} +
\big|y\big(\tau_{_{j(k)+1}}\big)- y\big({\tau_{\!_{j(k)+1}}}^{\!+}\big)
\big|
\Big)\qquad\qquad \forall~k\geq 1\,.
\tag 4.43
\endalign
$$
Thus, (4.39)-(4.41), (2.46) together, imply
$$
\align
\tau_2' -\tau_{\!_{j(2)}}
&\leq
\tau_2' -\tau_2
\\
\noalign{\smallskip}
&\leq \big|t_{\alpha_{i_1}}^{\,+}
\big(0,\,  
x^{\alpha_{i_1}}\big(\tau_2;\, 0,\, y(0^{\, +})\big)\big) \big|
\\
\noalign{\smallskip}
&\leq c_{10}\cdot \Big(
c_9\cdot \tv\{w\} \!+\!
\big|y(\tau_2)- y({\tau_2}^{\!+})
\big|
\Big)\,,
\tag 4.44
\endalign
$$
and
$$
\align
&\tau_k' -\tau_{\!_{j(k)}}\leq
\tau_k' -\tau_{\!_{j(k-1)+1}}
\\
\noalign{\smallskip}
&\quad \ \leq
\big|t_{\alpha_{i_{\!_{j(k)}}}}^{\,+}
\big(0,\,  
x^{\alpha_{i_{\!_{j(k)}}}}\big(\tau_{_{j(k)+1}};\, \tau_k',\, 
y({\tau'_k}^{\, +})\big)\big)
\big|+
\big|t_{\alpha_{i_{\!_{j(k-1)}}}}^{\,-}
\big(0,\,  y({\tau'_{k-1}}^{+})\big)
\big|
\\
\noalign{\smallskip}
&\quad\ \leq
c_{10}\cdot\Big(
c_9\cdot \tv\{w\} \!+\!
\big|y\big(\tau_{_{j(k-1)+1}}\big)- y\big({\tau_{\!_{j(k-1)+1}}}^{\!+}\big)
\big|
\Big)+
\big| t_{\alpha_{i_{\!_{j(k-1)}}}}^{\,-}
\big(0,\,  y({\tau'_{k-1}}^{+})\big)
\big|\qquad \forall~k>2\,.
\tag 4.45
\endalign
$$
On the other hand,  since \, $y(\cdot)$ \,
satisfies (1.3), using (4.31) we find
$$
\align
d\Big(y({\tau'_k}^{\, +}),\ 
B(\Om_{\alpha_{i_{\!_{j(k)}}}},\, \rho)\setminus 
D_{\alpha_{i_{\!_{j(k)}}}}\Big)
&\leq 
\big|y({\tau'_k}^{\, +})-y\big({\tau_{_{j(k)}}}^{\!+}\big)\big|+ 
\big|y\big({\tau_{_{j(k)}}}^{\!+}\big)-y\big(\tau_{_{j(k)}}\big)\big|
\\
\noalign{\smallskip}
&\leq 
M\big|\tau'_k-\tau_{_{j(k)}}\big|+\tv\{w\}
+ \big|y\big({\tau_{_{j(k)}}}^{\!+}\big)-y\big(\tau_{_{j(k)}}\big)\big|
\qquad \forall~k\,.
\tag 4.46
\endalign
$$
Therefore, 
proceeding by induction on \, $k,$ \, 
using (4.36) (thanks to (4.35), (4.38)), 
and relying on (4.37), (4.43)-(4.46), 
we obtain for every \, $k>1$ the  estimates
$$
\align
\tau_k' -\tau_{\!_{j(k)}}\leq
\tau_k' -\tau_{\!_{j(k-1)+1}}&\leq (2+ c_9)c_{10} \cdot (1 + M c_{10})^{k-2}
\cdot \tv\{w\}%
\\
\noalign{\smallskip}
&\leq \frac{c_{11}}{2(1+c_9)(1+2Mc_{10})} \cdot \tv\{w\}\,,
\tag 4.47
\endalign
$$
$$
\align
d\Big(y({\tau'_k}^{\, +}),\ 
B(\Om_{\alpha_{i_{\!_{j(k)}}}},\, \rho)\setminus 
D_{\alpha_{i_{\!_{j(k)}}}}\Big)
&\leq (2+ c_9) \cdot (1 + M c_{10})^{k-1}\cdot \tv\{w\}
\\
\noalign{\smallskip}
&\leq c_{11} \cdot \tv\{w\}<\delta_7\,
\tag 4.48
\endalign
$$
$$
\align
\big|t_{\alpha_{i_{\!_{j(k)}}}}^{\,-}
\big(0,\,  y({\tau'_k}^{\, +})\big)
\big|&\leq c_{10}\cdot d\Big(y({\tau'_k}^{\, +}),\ 
B(\Om_{\alpha_{i_{\!_{j(k)}}}},\, \rho)\setminus 
D_{\alpha_{i_{\!_{j(k)}}}}\Big)
\\
\noalign{\smallskip}
&\leq (2+ c_9)c_{10} \cdot (1 + M c_{10})^{k-1}
\cdot \tv\{w\}
\\
\noalign{\smallskip}
&\leq \frac{c_{11}}{2c_9(1+2Mc_{10})} \cdot \tv\{w\}\,.
\tag 4.49
\endalign
$$
\vs

\noindent
{\smbf 3.} \
Consider now the piecewise
smooth map \, $\widetilde y: [0,T]\mapsto\Bbb R^2$ \, 
defined by setting 
$$
\widetilde y(0)\doteq y(0),
\qquad\qquad\quad
\widetilde y(t) \doteq
x^{\alpha_{i_1}}\big(t;\ 0,\, y(0^{\, +}) \big)
\qquad\ \ t \in ]0,\,\tau'_2]\,,
\tag 4.50
$$
and
$$
\quad
\widetilde y(t) \doteq
x^{\alpha_{i_{\!_{j(k)}}}}\Big(t;\ \tau_k',\, 
x^{\alpha_{i_{\!_{j(k)}}}}
\big(t_{\alpha_{i_{\!_{j(k)}}}}^{\,-}
\big(0,\, 
y({\tau_k'}^{\, +})\big);\ 0,\,  y({\tau_k'}^{\, +}\big)
\big)\Big)
\qquad\quad t \in ]\tau_k',\,\tau_{k+1}']\,,
\quad\ 1<k\leq m'\,.
\tag 4.51
$$
By construction one has
$$
\aligned
\widetilde y({\tau_k'}^{\,+})&\in \partial D_{\alpha_{i_{\!_{j(k)}}}}
\qquad\quad \forall~1<k\leq m'\,,
\\
\noalign{\smallskip}
\widetilde y({\tau_{k+1}'})&\in \Omega\setminus 
\overset \,\circ \to D_{\alpha_{i_{\!_{j(k)}}}}
\qquad\ \forall~1\leq k<m'\,.
\endaligned
\tag 4.52
$$
Moreover, since
$\widetilde y(\cdot)$ is a classical solution of \, 
$\dot y = g_{\alpha_{i_{\!_{j(k)}}}}(y)$ \,
on  \, $]\tau'_k, \, \tau'_{k+1}[\,,$ 
and because $y(\cdot)$ satisfies~(1.3),
setting \, $t_{\alpha_{i_1}}^{\,-}
(0,\,y(0^{\, +}))\doteq 0,$ \, 
and using  (4.31), (4.34),  (4.47), (4.49), 
we derive, for all \, $k,$ \, the estimates
$$
\align
\big|\widetilde y(t)-y(t)\big|&\leq
\big|\widetilde y(t)-y({\tau'_k}^{\, +})\big|+
\big|y({\tau'_k}^{\, +})-y(t)\big|
\\
\noalign{\smallskip}
&=\big|\widetilde y(t)-\widetilde y(\tau_k'\!-\!
t_{\alpha_{i_{\!_{j(k)}}}}^{\,-}
\!\big(0,\,  y({\tau'_k}^{\, +})\big))\big|+
\big|y({\tau'_k}^{\, +})-y(t)\big|
\\
\noalign{\smallskip}
&\leq 2M\cdot \big| t_{\alpha_{i_{\!_{j(k)}}}}^{\,-}
\big(0,\,  y({\tau'_k}^{\, +})\big)
\big| + \tv\{w\}
\\
\noalign{\smallskip}
&\leq c_{11} \cdot \tv\{w\}
\\
&\qquad\qquad\qquad \forall~t \in \big]\tau_k',\, \tau_k'\!-\!
 t_{\alpha_{i_{\!_{j(k)}}}}^{\,-}
\!\big(0,\,  y({\tau'_k}^{\, +})\big)\!\big]\,,
\tag 4.53
\endalign
$$
$$
\align
\big|\widetilde y(t)-y(t)\big|&=\Big|
x^{\alpha_{i_{\!_{j(k)}}}}\Big(t;\ \tau_k'\!-\!
 t_{\alpha_{i_{\!_{j(k)}}}}^{\,-}
\!\big(0,\, y({\tau_k'}^{\, +})\big),\,
y({\tau_k'}^{\, +})\Big)-
z^{\alpha_{i_{\!_{j(k)}}}}\big(t;\ w,\,\tau_k',\,
y({\tau_k'}^{\, +}\big)\big)
\Big|
\\
\noalign{\smallskip}
&\leq c_9\cdot \Big(\big|t_{\alpha_{i_{\!_{j(k)}}}}^{\,-}
\big(0,\,  y({\tau'_k}^{\, +})\big)
\big| + \tv\{w\}\Big)
\\
\noalign{\smallskip}
&\leq c_{11} \cdot \tv\{w\}
\\
&\qquad\qquad\qquad \forall~t \in \big]\tau_k'\!-\!
t_{\alpha_{i_{\!_{j(k)}}}}^{\,-}
\!\big(0,\,  y({\tau'_k}^{\, +})\big),\,
\tau_{\!_{j(k)+1}}\big]\,,
\tag 4.54
\endalign
$$
$$
\align
\big|\widetilde y(t)-y(t)\big|&\leq
\big|\widetilde y(t)-\widetilde y(\tau_{\!_{j(k)+1}})\big|+
\big|\widetilde y(\tau_{\!_{j(k)+1}})-y(\tau_{\!_{j(k)+1}})\big|+
\big|y({\tau_{\!_{j(k)+1}}}^{\!+})-y(\tau_{\!_{j(k)+1}})\big|+
\\
\noalign{\smallskip}
&\qquad\quad + \big|y(t)-y({\tau_{\!_{j(k)+1}}}^{\!+})\big|
\\
\noalign{\smallskip}
&\leq 2M\cdot \big|\tau_{k+1}'-\tau_{\!_{j(k)+1}}\big|+ c_9\cdot
\big|t_{\alpha_{i_{\!_{j(k)}}}}^{\,-}
\big(0,\,  y({\tau'_k}^{\, +})\big)\big|
+(2+ c_9)\cdot \tv\{w\}
\\
\noalign{\smallskip}
&\leq (2M+c_9)\cdot \big|\tau_{k+1}'-\tau_{\!_{j(k)+1}}\big|
+(2+ c_9)\cdot \tv\{w\}
\\
\noalign{\smallskip}
&\leq c_{11} \cdot \tv\{w\}
\\
&\qquad\qquad\qquad 
\forall~t \in \big]\tau_{\!_{j(k)+1}},\, \tau_{k+1}'\big]\,,
\tag 4.55
\endalign
$$
which, together, and thanks to (4.38), yield
$$
\big\|\widetilde y-y\big\|_{\L^\infty([0,T])}
\leq c_{11}\cdot \tv\{w\}<\rho\,.
\tag 4.56
$$
\v

\noindent
{\smbf 4.} \
Because of (4.56), the choice of the constants 
\, $C^{\text{\sl v\i}}\geq 2c_{11},\, C^{\text{\sl v\i\i}},\,
\delta^{\text{\sl v\i}}$ \, according
with Lemma~4.1 guarantees that
$$
\text{\rm meas}\,
\left(\bigcup_k
\Big\{
t\in ]\tau_k',\,\tau_{k+1}']~:~
\widetilde y(t)\notin D_{\alpha_{i_{\!_{j(k)}}}}
\Big\}\right)
< C^{\text{\sl v\i\i}} \cdot \tv\{w\}.
\tag 4.57
$$
Therefore, relying on (4.52), (4.57),
the  inward-pointing
condition (1.15) together with
the transversality condition {\bf (C)} imply that,
for every \, $k=1,\dots,m',$ \, there exists
a partition \, $\tau_{k}'=\widetilde t_{k,1}<\widetilde t_{k,1}
<\cdots < \widetilde t_{k,\widetilde q_k}=\tau_{k+1}'$ \, of \,
$[\tau'_k, \tau'_{k+1}],$ \, so that 
$$
\gathered
\widetilde y(t) \in D_{\alpha_{i_{\!_{j(k)}}}} \qquad
\quad \forall ~t \in ]\,\widetilde t_{k,\ell},\, \widetilde t_{k,\ell+1}\,[
\qquad\quad \text{for \ all \ odd}\quad 
\ell,
\\
\noalign{\medskip}
\sum_{\ell \  \text{even}} (\,\widetilde t_{k,\ell+1}-\widetilde t_{k,\ell})
\leq C^{\text{\sl v\i\i}} \cdot \tv\{w\}\,,
\endgathered
\tag 4.58
$$
and with the property that the points \, 
$$
\widetilde y(\,\widetilde t_{k,1}^{\ \ \,+}),\quad t_{k,1} \neq 0,
\qquad\qquad \widetilde y(\,\widetilde t_{k,\ell}), \quad 1 < \ell 
\leq 2 \lfloor \widetilde q_k/2\rfloor,
$$ 
lie on different edges of the domain \, $D_{\alpha_{i_{\!_{j(k)}}}}$ \,
($\lfloor a \rfloor$ denoting the integer part of
$a$). \,
For each \, $k=1,\dots,m',$ \, consider  the sequence of points
\, $t_{k,1}< t_{k,2}<\cdots <t_{k,q_k},$ \ 
$q_k\doteq\lfloor \widetilde q_k/2\rfloor+1,$ \, 
recursively
defined by setting \, $t_{1,1}\doteq 0,$ \, and 
$$
\gathered
t_{k,1}
\doteq t_{k-1,\,q_{k-1}}
\qquad\quad 1<k\leq m'\,,
\\
\noalign{\medskip}
t_{k,\ell+1}
\doteq t_{k,\ell}+ \widetilde t_{k,2\ell}-
\widetilde t_{k,2\ell-1}
\qquad\quad 1\leq \ell < q_k\,,
\qquad 1\leq k\leq m'\,.
\endgathered
\tag 4.59
$$
Then, letting
\, $x^{g}\big(t;\, t_0,\,x_0)$ \,
denote a solution of (1.11) starting 
from \, $x_0$ \, at time \, $t_0,$ \,
we define the map
\, $y^\di: [0,T]\mapsto\Bbb R^2$ \, as follows: 
$y^\di(0)\doteq \widetilde y(0),$ and
$$
y^\di(t)\doteq
\cases
\widetilde y\big(\,t+\widetilde t_{k,2\ell-1}-t_{k,\ell}\,\big)
\quad\quad 
&\forall~t\in ]\,t_{k,\ell},\,t_{k,\ell+1}\,]\,,
\qquad\quad 1\leq\ell < q_k\,,
\qquad 1\leq k \leq m'\,,
\\
\noalign{\medskip}
x^{g}\big(t;\, t_{m',q_{m'}},\,y^\di(t_{m',q_{m'}})\big)
\quad\quad 
&\forall~t\in ]\,t_{m',q_{m'}},\,T\,]\,.
\endcases
\tag 4.60
$$
Notice that, by construction,
and by the properties of the
solutions of a patchy system 
(recalled in Section~1),  the map \, $t \mapsto y^\di(t)$ \,
enjoys the properties $a')$-$b')$ stated in Lemma~2.2.
Moreover, since one has
$$
\widetilde t_{k,2\ell-1}-t_{k,\ell}
\leq \sum_{p \  \text{even}} 
\widetilde t_{k,p+1}-\widetilde t_{k,p}
= T-t_{m',q_{m'}}
\qquad\quad \forall~l,\ \ \forall~k\,,
$$
and because $\widetilde y(\cdot)$ is a solution of \, 
$\dot y = g_{\alpha_{i_{\!_{j(k)}}}}(y)$ \,
on  \, $]\tau'_k, \, \tau'_{k+1}[\,,$ 
using (4.31), (4.56), (4.58), we derive
$$
\align
\big|y^\di(t)-\widetilde y(t)\big|
&\leq \!\!\sum_{1< k\leq m'} 
\big|\widetilde y({\tau'_k}^{\, +})-\widetilde y(\tau'_k)\big|+
M \cdot \max_{k,\ell} \big|\,\widetilde t_{k,2\ell-1}-t_{k,\ell} \big|
\\
&\leq \!\!\!\sum_{1< k\leq m'} \!\!
\Big(\big|\widetilde y({\tau'_k}^{\, +})-y({\tau'_k}^{\, +})\big|\!+\!
\big|y({\tau'_k}^{\, +})- y(\tau'_k)\big|\!+\!
\big|\widetilde y(\tau'_k)- y(\tau'_k)\big|
\Big)\!\!+\!
M C^{\text{\sl v\i\i}}\cdot \tv\{w\}
\\
\noalign{\medskip}
&\leq  \big(1+2c_{11}+ MC^{\text{\sl v\i\i}}\,\big) \cdot \tv\{w\}
\qquad \qquad \forall~t\in[0,\,t_{m',q_{m'}}]\,,
\tag 4.61
\endalign
$$
%

\noindent
and
$$
\align
\big|y^\di(t)-\widetilde y(t)\big|
&\leq \big|y^\di(t)-y^\di(t_{m',q_{m'}})\big|+
\big|y^\di(t_{m',q_{m'}})-\widetilde y (t_{m',q_{m'}})\big|+
\big|\widetilde y(t)-\widetilde y(t_{m',q_{m'}})\big| 
\\
\noalign{\medskip}
&\leq \sum_{1< k\leq m'} 
\big|\widetilde y({\tau'_k}^{\, +})-\widetilde y(\tau'_k)\big|+
2M \cdot \big|T-t_{m',q_{m'}}\big| +
\big|y^\di(t_{m',q_{m'}})-\widetilde y (t_{m',q_{m'}})\big|
\\
\noalign{\smallskip}
&\leq  \big(2+4c_{11}+ 3MC^{\text{\sl v\i\i}}\,\big) \cdot \tv\{w\}
\qquad \quad \forall~t\in[t_{m',q_{m'}},\,T]\,.
\tag 4.62
\endalign
$$
\v
\noindent
Hence, (4.61)-(4.62) together with (4.56), show that \, $y^\di(\cdot)$ \,
satisfies the estimates (2.7)-(2.8) of  property~$c')$,
taking 
$$
C''> \big(5+10 c_{11}+ 6MC^{\text{\sl v\i\i}}\,\big)\,,
\tag 4.63
$$
which completes the proof of Lemma~2.2.
\fine
\vs

\noindent
{\bf Proof of Lemma~2.3.}  \ 
\v

\noindent 
{\smbf 1.} \ 
As in the proof of Lemma~4.1, call \, 
$E^1, \dots ,E^{^{p_{\alpha_o}}},$ \, 
and \, $r^1, \dots ,r^{^{p_{\alpha_o}}},$ \, 
respectively,
the edges of  the polygonal
domain \, $D_{\alpha_o}$ \,
that form the boundary~$\partial D_{\alpha_o},$ \, 
and the corresponding lines in which the edges are contained.
Let  \, $\bn^{{\boldsymbol \ell}}$ \,  
be the  normal to~$E^\ell$ \, pointing towards 
the interior of \, $D_{\alpha_o},$ \,
and call 
\, $\Cal I_{\alpha_o}$ \, and \, $\Cal O_{\alpha_o},$ \, respectively, 
the set of incoming and of outgoing indices defined as in (4.12),
so that
$$
\alignedat 4
\big \langle
g_{\alpha_o}(x),\ 
\bn^{{\boldsymbol \ell}}
\big\rangle
&>0
\qquad\quad &\forall~x \in 
\overline D_{\alpha_o}\cap r^\ell\,,
\qquad\quad \forall~\ell \in \Cal I_{\alpha_o}\,,
\\
\noalign{\medskip}
\big \langle
g_{\alpha_o}(x),\ 
\bn^{{\boldsymbol \ell}}
\big\rangle
&<0
\qquad\quad &\forall~x \in 
\overline D_{\alpha_o}\cap r^\ell\,,
\qquad\quad \forall~\ell \in \Cal O_{\alpha_o} \,.
\endalignedat
\tag 4.64
$$
Let \, $V_{\alpha_o}$ \, be the set of vertices of the domain \,
$D_{\alpha_o},$ \, and denote
$$
\Cal T \doteq \Big\{
\gamma_j : [0,\, \widehat \tau^j\,] \mapsto\Bbb R^2~;~
j=1, \dots, m
\Big\}
\tag 4.65
$$
the set of maximal trajectories of \, $\dot y = g_{\alpha_o}(y)$ \, 
that go through some vertex in \, $V_{\alpha_o},$ 
and whose graph is contained in \, $\overline D_{\alpha_o}.$
Thus, we have
$$
\gathered
\text{Im}(\gamma_j)\cap V_{\alpha_o}\neq \emptyset\,,
\qquad\quad 
\gamma_j(0),\, \gamma_j(\widehat \tau^j)\in\partial D_{\alpha_o}\,,
\\
\noalign{\smallskip}
\gamma_j(t)=x^{\alpha_o}\big(t;\ 0,\, \gamma_j(0)\big)
\in \overline D_{\alpha_o}
\qquad\quad\forall~t\in [0,\, \widehat \tau^j\,]\,.
\endgathered
\tag 4.66
$$
In connection with every trajectory \, $\gamma_j,$ \,
there will be a partition \, $\widehat t_{j,1}=0<\widehat t_{j,2}<
\cdots < \widehat t_{j,k_j}=\widehat \tau^j$ \,
of \, $[0,\, \widehat \tau^j\,]$
such that
$$
\gathered
\gamma_j(\,\widehat t_{j,h})\in \partial D_{\alpha_o}
\qquad\quad \forall~1 \leq h \leq k_j\,,
\\
\noalign{\medskip}
\gamma_j(t)\in \overset \,\circ \to D_{\alpha_o}
\qquad\ \forall~t\in ]\,\widehat t_{j,h},\, \widehat t_{j,h+1}[\,,
\qquad 1 \leq h <  k_j\,.
\endgathered
\tag 4.67
$$
Then, set
$$
\gamma_{_{j,h}}\doteq 
\gamma_j \restriction_{[\,\widehat t_{j,h},\, \widehat t_{j,k_j}]},
$$ 
and let
$$
\aligned
{\Cal L}^{\Cal I}_{j,h}
&\doteq \big\{\ell \in {1, \dots , p_{\alpha_o}}\ ;\ \
\gamma_j(\,\widehat t_{j,h})\in E^\ell
\qquad\text{and}\qquad \ell \in \Cal I_{\alpha_o}\big\}\,,
\\
\noalign{\medskip}
{\Cal L}^{\Cal O}_{j,h}
&\doteq \big\{\ell \in {1, \dots , p_{\alpha_o}}\ ;\ \
\gamma_j(\,\widehat t_{j,h})\in E^\ell
\qquad\text{and}\qquad \ell \in \Cal O_{\alpha_o}\big\}\,,
\endaligned
\tag 4.68
$$
denote the incoming and outgoing indices
of the edges of \, $D_{\alpha_o}$ \,
that pass through \, $\gamma_j(\widehat t_{j,h}).$ \,
Moreover, set
$$
{\Cal L}_{j,h}
\doteq {\Cal L}^{\Cal I}_{j,h} \cup {\Cal L}^{\Cal O}_{j,h}\,.
\tag 4.69
$$
Since we are assuming that \, $g_{\alpha_o}$ \, is a smooth
vector field defined on a neighborhhod of \, $\overline D_{\alpha_o}$ \,
that satisfies the transversality condition {\bf (C)} and 
the inward-pointing
condition~(1.15),
there will be some constant \, $\delta_8>0$ \, so that one has
$$
\forall~x_0\in \overline D_{\alpha_o}
\cap B\big(\gamma_j(\,\widehat t_{j,h})
,\, \delta_8\big)\,,\quad
\forall~h'>h\,,\quad
\forall~\ell\in {\Cal L}^{\Cal O}_{j,h'}
\qquad\exists !~~~~t>0
\quad : \quad 
x^{\alpha_o}\big(t;\ 0,\, x_0
\big)\in r^\ell\,.
\tag 4.70
$$
%
Thus, for such \, $x_0, \, \ell,$ \,
define
$$
%
t^\ell(x_0)\doteq
t>0\quad\ \text{s.t.}\quad x^{\alpha_o}(t;\ 0,\, x_0)
\in r^{\ell}
\,.
\tag 4.71
$$
The regularity properties of the flow map of \, $g_{\alpha_o}$ \,
together with the transversality condition {\bf (C)} and 
the inward-pointing
condition~(1.15)
guarantee that there exists some constant \, $c_{12}>1$ \, such that
$$
\big|t^\ell(x_0)-t^\ell(y_0)\big|\leq c_{12} \cdot 
\big|x_0-y_0\big|\qquad\quad \forall~x_0,\, y_0\,.
\tag 4.72
$$
\vs

\noindent
{\smbf 2.} \
We shall construct now an increasing tube 
\, $\Gamma(\gamma_j, \, \lambda)$ \,
of size \, $\lambda$ \,
around each  trajectory \,
$\gamma_j\in \Cal T,$ \, which is positively invariant w.r.t.
left-continuous, piecewise smooth
function  having the properties~$a'')$-$b'').$
Take
\, $0<\lambda<\delta_8,$ \, and
let
$$
\gathered
F(\gamma_{_{j,1}}, \, \lambda)\doteq 
B\Big(\gamma_j(\,\widehat t_{j,1}),\, \lambda\Big)\cap 
\Bigg(\bigcup_{\ell \in {\Cal L}
_{j,1}} E^\ell\Bigg) \,,
\\
\noalign{\medskip}
G(\gamma_{_{j,1}}, \, \lambda)\doteq
\Big\{
x^{\alpha_o}(s;\ 0,\, x_0)\quad;\ \
x_0\in F(\gamma_{_{j,1}}, \, \lambda),
\quad\ 0\leq s \leq t^{\ell'}(x_0),
\quad \ell' \in {\Cal L}^{\Cal O}_{j,k_j}
\Big\}\cap D_{\alpha_o}\,.
\endgathered
\tag 4.73
$$
Then, proceeding by induction on $h>1,$ and relying on 
(4.34), (4.70), (4.72),
one can show that 
there exists  constants \, $c_{13}>1,$ \, $\overline \lambda<\delta_8/c_{13},$ \,
such that, 
letting 
%
$$
\gathered
F(\gamma_{_{j,h}}, \, \lambda)\doteq 
B\Big(G(\gamma_{_{j,h-1}}, \, \lambda),\, \lambda\Big)\cap 
\Bigg(\bigcup_{\ell \in {\Cal L}
_{j,h}} E^\ell\Bigg) \,,
\\
\noalign{\medskip}
G(\gamma_{_{j,h}}, \, \lambda)\doteq
\Big\{
x^{\alpha_o}(s;\, 0, x_0)\ \ ;\ \,
x_0\in F(\gamma_{_{j,h}}, \, \lambda),
\quad\ 0\leq s \leq t^{\ell'}(x_0),
\quad \ell' \in {\Cal L}^{\Cal O}_{j,k_j}
\Big\}\cap D_{\alpha_o}\,,
\endgathered
\quad \lambda\leq \overline \lambda\,,
\tag 4.74
$$
for \, $1< h < k_j,$ \, 
one has
$$
\big|\gamma_j(\,\widehat t_{j,h})-x_0\big|\leq c_{13} \cdot \lambda
\qquad\quad \forall~x_0\in F(\gamma_{_{j,h}}, \, \lambda),
\quad 1\leq  h < k_j, \quad \lambda \leq \overline\lambda\,.
\tag 4.75
$$
%
Moreover, we may choose \, $\overline \lambda$ \,
so that, setting
$$
\Gamma(\gamma_j, \, \lambda)\doteq
\bigcup_{h=1}^{k_j-1}
G(\gamma_{_{j,h}}, \, \lambda)\,,
\qquad\quad 1\leq j\leq m\,,
\tag 4.76
$$
there holds
$$
\Gamma(\gamma_j, \, \lambda)\cap \Gamma(\gamma_i, \, \lambda)=
\emptyset\,,\qquad\quad \forall~1\leq i,j \leq m, \quad\ i\neq j,
\qquad\quad \forall~\lambda \leq \overline \lambda\,.
\tag 4.77
$$
Let \, 
$R^1(\lambda), \dots, R^r(\lambda),$ \,
denote the connected components of \, 
$\overline 
D_{\alpha_o}\!\setminus 
\cup_{j=1}^m\Gamma(\gamma_j, \, \lambda),$ \,
and set
$$
\rho
(\lambda)\doteq \min
\Big\{ d\big(R^s(\lambda),\,
E^{\ell}\big)\ \ ;\ \, 1 \leq s \leq r,\ \ \
1 \leq \ell \leq p_{\alpha_o}, \ \ \ 
R^s(\lambda) \cap E^\ell = \emptyset 
\Big\}.
\tag 4.78
$$
Observe that, 
by the transversality condition {\bf (C)} and because
of the inward-pointing
condition~(1.15),  one has
$$
\inf
\left\{\,
\frac{\rho
(\lambda)}{\lambda}~:~ 0<\lambda \leq \overline \lambda
\,\right\}>0.
\tag 4.79
$$
Therefore, 
there exist 
constants \, $c_{14}>1, \ 0<\delta_9<\overline \lambda/c_{14}, $ \,
so that
$$
\rho
(c_{14} \cdot \delta) > 2 \delta
\qquad\quad \forall~0<\delta\leq \delta_9.
\tag 4.80
$$
\vs

\noindent
{\smbf 3.} \
Consider now a left-continuous, piecewise smooth
function \, $y^\bem :~]\tau_0, \tau_1]\mapsto\Bbb R^2$ \, 
having the properties $a'')$-$b'')$ with
$$
\overline \delta \leq \delta_9\,.
\tag 4.81
$$
To fix the ideas, assume that Case~$a''$-{\it 3})
occurs which, in particular, implies
$$
y^\bem(\tau_1)\in \partial D_{\alpha_o}\,.
\tag 4.82
$$
The other two cases can be treated in an entirely similar manner.
Since, by construction, the tubes \, $G(\gamma_{j,h}, \, \lambda)$ \,
defined in (4.73)-(4.74) are 
invariant subsets of~$\overline D_{\alpha_o}$ \,
for the flow map
of \, $g_{\alpha_o},$ \,   if
$$
y^\bem(\tau_0^{\, +})\in D_{\alpha_o}\setminus
\bigcup_{j=1}^m
\Gamma(\gamma_j, \, c_{14} \cdot \Delta(y^\bem)\big)\,,
\tag 4.83
$$
then, by property $a'')$ it follows that
$$
y^\bem(t)\in \overline D_{\alpha_o}\setminus
\bigcup_{j=1}^m
\Gamma(\gamma_j, \, c_{14} \cdot \Delta(y^\bem)\big)
\qquad\quad\forall~t\in ~]t'_1,\, t'_2]\,.
\tag 4.84
$$
We claim that (4.84) implies \, $q_o=2.$ \,
Indeed, if \, $R^s(c_{14} \cdot \Delta(y^\bem))$ \, 
denotes the connected component of \, 
$\overline D_{\alpha_o}\!\setminus 
\big(\cup_{j=1}^m\Gamma(\gamma_j, \, c_{14} \cdot \Delta(y^\bem))\big)$ \,
that contains \, $\text{Im}(y^\bem\!\restriction_{]t'_1,\, t'_2]}),$ \,
and we assume by contradiction that \, $q_o>2,$ \, relying on~(4.77)
and on property $a''$-{\it 3})
we deduce that
\, $y^\bem({t'_2}^{\, +})$ \, lies on some edge \, $E^{\ell}$ \,
such that 
$R^s(c_{14} \cdot \Delta(y^\bem))
\cap E^{\ell}= \emptyset.$ \,
But then, using the inequality~(4.80), and because of the
definition (4.78) of the quantity \, $\rho
,$ \,
we would obtain
$$
\align
\big|y^\bem(t'_2)-y^\bem({t'_2}^{\, +})\big|
&\geq d\big(R^s(c_{14} \cdot \Delta(y^\bem)),\, 
E^{\ell}\,\big)
\\
&\geq \rho
\big(c_{14} \cdot \Delta(y^\bem)\big)
\\
&>2 \Delta(y^\bem)
\endalign
$$
which yields a contradiction by the definition of 
\, $\Delta(y^\bem)$ \, at (2.10). 
Therefore, it must be \, $q_0=2,$ and hence
we have 
$$
y^\bem(t)=
x^{\alpha_o}\big(t;\, \tau_0, \, y^\bem(\tau_0^{\, +})\big)
\qquad\quad \forall~t \in ]\tau_0,\, \tau_1]\,,
$$
which, together with (2.9), (4.82), clearly shows that
\, $Q_{\alpha_o}\doteq y^\bem(\tau_0^{\, +}),$ \, 
\, $\sigma_{\alpha_o}\doteq \tau_1$ \, enjoy pro-\linebreak
perties~$c'')$-$e'')$ whenever (4.83) holds.
\vs

\noindent
{\smbf 4.} \
Assume now that (4.83) is not verified
i.e. that, for some \, $1\leq j \leq m,$ \,
$1 \leq h' \leq k_j,$ \, there holds
$$
y^\bem(\tau_0^{\, +})\in 
G\big(\gamma_{_{j,h'}}, \, c_{14} \cdot \Delta(y^\bem)\big)\,.
\tag 4.85
$$
Proceeding by induction on \, $h' \leq h < k_j,$ \,
and relying on (4.70), (4.77) and on property $a'')$,
one then easily derives that 
$$
y^\bem(t)\in \bigcup_{h=h'}^{k_j-1}
G\big(\gamma_{_{j,h}}, \, c_{14} \cdot \Delta(y^\bem)\big)
\subset \Gamma\big(\gamma_j, \, c_{14} \cdot \Delta(y^\bem)\big)
\qquad\quad \forall~t\in~]\tau_0, \tau_1]\,.
\tag 4.86
$$
%
In connection with the partition 
\, $t'_1=\tau_0<t'_2<\cdots<t'_{q_o}=\tau_1$ \, of \, $[\tau_0,\tau_1],$  
\, induced by ~$y^\bem(\cdot)$ \,
according with property $a'')$,
define the triplet
of indices \, $h(\ell),\, p^-(\ell),\, p^+(\ell),$ \, $1<\ell \leq q_o$ \,
so that 
\, $E^{p^-(\ell)},$ \,  \, $E^{p^+(\ell)}$ \, denote
the edges of \, $D_{\alpha_o}$ \, which cross the trajectory
\, $\gamma_j$ \, in \, $\gamma_j(\,\widehat t_{j,h(\ell)}),$ \,
$h(\ell)>h',$ \,
and contain, respectively, the
point  \, $y^\bem(t'_\ell),$ \,
and the point \, $y^\bem({t'_\ell}^{\, +}),$ \, i.e.
such that
$$
\alignedat 4
y^\bem(t'_\ell)&\in E^{p^-(\ell)}, \qquad 
{p^-(\ell)}&\in {\Cal L}^{\Cal O}_{j,h(l)}\,,
\qquad\quad \forall~1<\ell \leq q_0\,,
\\
\noalign{\smallskip}
y^\bem({t'_\ell}^{\, +})&\in E^{p^+(\ell)}, \qquad 
{p^+(\ell)}&\in {\Cal L}^{\Cal I}_{j,h(l)}\,
\qquad\quad \forall~1<\ell<q_0\,.
\endaligned
\tag 4.87
$$
%
%
Then, set
$$
Q_{\alpha_o}\doteq x^{\alpha_o}\big(
t_1'-t_2';\ 0,\, 
\gamma_j(\,\widehat t_{j,h(2)})\big)\,,
\tag 4.88
$$
and observe that, by property $a'')$,
one has 
$$
y^\bem(\tau_0^{\, +})
=x^{\alpha_o}\big(t_1'-t_2';\ 0,\, y^\bem(t_2')\big)\,.
\tag 4.89
$$
%
%
%
%
Thus,
using (4.34),  (4.88)-(4.89), we obtain
%
$$
\big|Q_{\alpha_o}-y^\bem(\tau_0^{\, +})\big|
\leq c_9\cdot
\big|\gamma_j(\,\widehat t_{j,h(2)})-
y^\bem(t_2')\big|\,.
\tag 4.90
$$
On the other hand, since one can show with
an inductive argument that (4.86)-(4.87) imply
$$
\aligned
y^\bem(t'_\ell)&\in 
F\big(\gamma_{_{j,h}}, \, c_{14} \cdot \Delta(y^\bem)\big)\,,
\qquad \forall~1<\ell\leq q_0\,,
\\
\noalign{\smallskip}
y^\bem({t'_\ell}^{\, +}) &\in 
F\big(\gamma_{_{j,h}}, \, c_{14} \cdot \Delta(y^\bem)\big)\,,
\qquad \forall~1<\ell<q_0\,,
\endaligned
$$
relying on  (4.75) we derive
$$
\aligned
\big|y^\bem(t'_\ell)-\gamma_j(\,\widehat t_{j,h(\ell)})
\big| &\leq c_{13}\cdot c_{14} \cdot \Delta(y^\bem)\,,
\qquad\ \forall~1<\ell\leq q_0\,,
\\
\noalign{\medskip}
\big|y^\bem({t'_\ell}^{\, +})-\gamma_j(\,\widehat t_{j,h(\ell)})
\big| &\leq c_{13}\cdot c_{14} \cdot \Delta(y^\bem)\,,
\qquad\ \forall~1<\ell<q_0\,.
\endaligned
\tag 4.91
$$
%
In turn, the first estimate in (4.91) for $\ell=2,$
together with (4.90), yields
$$
\big|Q_{\alpha_o}-y^\bem(\tau_0^{\, +})\big|
\leq c_9\cdot c_{13}\cdot c_{14} \cdot \Delta(y^\bem)\,.
\tag 4.92
$$
Moreover,  observing that by the definitions  (4.71), (4.87),
one has
$$
\aligned
t'_{\ell+1}-t'_\ell &= t^{p^-(\ell+1)}\big(y^\bem(t'_\ell)\big)\,,
\\
\noalign{\medskip}
\widehat t_{j,h(\ell+1)}-\widehat t_{j,h(\ell)} &=
t^{p^-(\ell+1)}\big(\gamma_j(\,\widehat t_{j,h(\ell)})\big)\,,
\endaligned
\qquad\quad \forall~1<\ell<q_0\,,
$$
thanks to (4.72) we get
$$
\big|
(t'_{\ell+1}-t'_\ell)-(\widehat t_{j,h(\ell+1)}-\widehat t_{j,h(\ell)})
\big|\leq c_{12}\cdot c_{13}\cdot c_{14} \cdot \Delta(y^\bem)\,,
\qquad\quad \forall~1<\ell<q_0\,.
\tag 4.93
$$
Therefore, since by the definition (4.88) we have
$$
Q_{\alpha_o}=\gamma_j\big(\,\widehat t_{j,h(2)}+t_1'-t_2'\big)\,,
$$
which, in particular, implies
$$
\gathered
x^{\alpha_o}\big(t;\, \tau_0, \, Q_{\alpha_o}\big)=
\gamma_j\big(\,\widehat t_{j,h(2)}-t_2'+t\big)\,,
\\
\noalign{\medskip}
x^{\alpha_o}\big(\,\widehat t_{j,h(\ell)}-\widehat t_{j,h(2)}
+t_2';\ \tau_0, \, Q_{\alpha_o}\big)=
\gamma_j\big(\,\widehat t_{j,h(\ell)}\big)
\qquad\forall~\ell\geq 2,
\endgathered
\tag 4.94
$$
and because by property $a'')$
one has 
$$
y^\bem(t)
=x^{\alpha_o}\big(t;\ t'_\ell,\, y^\bem({t'_\ell}^{\, +})\big)
\qquad\quad\forall~t\in~]t'_\ell,\,t'_{\ell+1}]\,,
\qquad 1\leq \ell<q_0\,,
\tag 4.95
$$
relying on (4.91)-(4.93), and using
(4.34), (4.94)-(4.95), we derive
$$
\align
&\big|y^\bem(t)-
x^{\alpha_o}\big(t;\, \tau_0, \, Q_{\alpha_o}\big) \big| \leq
\\
\noalign{\smallskip}
&\quad \leq
\big|x^{\alpha_o}\big(t;\  t'_\ell,\, y^\bem({t'_\ell}^{\, +})\big)-
x^{\alpha_o}\big(t;\  t'_\ell,\, 
\gamma_j\big(\,\widehat t_{j,h(\ell)}\big)\big)
\big|+
\big|x^{\alpha_o}\big(t;\  t'_\ell,\, 
\gamma_j\big(\,\widehat t_{j,h(\ell)}\big)\big)-
x^{\alpha_o}\big(t;\, t_1', \, Q_{\alpha_o}\big)\big|
\\
\noalign{\smallskip}
&\quad = 
\big|x^{\alpha_o}\big(t;\  t'_\ell,\, y^\bem({t'_\ell}^{\, +})\big)-
x^{\alpha_o}\big(t;\  t'_\ell,\, 
\gamma_j\big(\,\widehat t_{j,h(\ell)}\big)\big)
\big|+
\\
\noalign{\smallskip}
&\qquad \quad +
\big|x^{\alpha_o}\big(t+\!(\,\widehat t_{j,h(\ell)}\!-\widehat t_{j,h(2)})\!-
\!(t'_\ell-\!t'_2);\  t'_1,\, 
Q_{\alpha_o}\big)-
x^{\alpha_o}\big(t;\, t_1', \, Q_{\alpha_o}\big)\big|
\\
\noalign{\smallskip}
&\quad \ \leq c_9\cdot \Big(
\big|y^\bem(t'_\ell)-\gamma_j(\,\widehat t_{j,h(\ell)})
\big|+
\big|(\,\widehat t_{j,h(\ell)}\!-\widehat t_{j,h(2)})\!-
\!(t'_\ell-\!t'_2)\big|
\Big)
\\
\noalign{\smallskip}
&\quad \ \leq c_9\cdot c_{13}\cdot c_{14}\cdot\big(1+ q_0\cdot c_{12}\big)
\cdot \Delta(y^\bem)\,,
\\
\noalign{\medskip}
&\hskip 1.5in
\forall~t\in ~]t_\ell',\ 
\min\{t'_{\ell+1},\ \widehat t_{j,h(q_0)}\!-\widehat t_{j,h(2)}\!+t_2'\}],
\qquad 1\leq \ell<q_0\,,
\tag 4.96
\endalign
$$
and
$$
\align
\big|\big(\,\widehat t_{j,h(q_0)}\!-
\widehat t_{j,h(2)}\!+t_2'\big)-\tau_1\big|
&=
\big|(\,\widehat t_{j,h(q_0)}\!-\widehat t_{j,h(2)})\!-
(t_{q_0}'-t_2')\big|
\\
\noalign{\smallskip}
&\leq q_0\cdot c_{12} \cdot c_{13}\cdot c_{14}\cdot \Delta(y^\bem)\,.
\tag 4.97
\endalign
$$
Moreover,  
if we set \, $\sigma_{\alpha_o}\doteq 
\widehat t_{j,h(q_0)}\!-\widehat t_{j,h(2)}\!+t_2',$ \,
thanks to (4.67), (4.94) we obtain
$$
x^{\alpha_o}\big(\sigma_{\alpha_o};\, \tau_0, \, Q_{\alpha_o}\big)=
\gamma_j\big(\,\widehat t_{j,h(q_0)}\big)
\in \partial D_{\alpha_o}\,.
\tag 4.98
$$
Hence, 
(4.96)-(4.98) together, show that \, $Q_{\alpha_o}$ \,
defined as in (4.88)
and \, $\sigma_{\alpha_o}=
\widehat t_{j,h(q_0)}\!-\widehat t_{j,h(2)}\!+t_2'$ \,
enjoy properties~$c'')$-$e''),$ 
taking \, $\overline C>
c_9\cdot c_{13}\cdot c_{14}\cdot\big(1+ q_0\cdot c_{12}\big),$ \,
in the case where (4.85) is verified,
which completes the proof of Lemma~2.3.
\fine

%
\parindent 60pt

\vsk
\centerline{\medbf References.}
\vs
\noindent\item{\ABpv\ } {\smc F. Ancona, A. Bressan},
Patchy vector fields and asymptotic stabilization,
{\it ESAIM - Control, Optim. and Calc. of Var.,} 
Vol. {\bf 4}, (1999), pp. 445-471.
\medskip
\noindent\item{\ABfs\ } {\smc F. Ancona, A. Bressan},
Flow Stability of Patchy vector fields and 
Robust Feedback Stabilization, (2001),
submitted.
\medskip
%
%
%
%
\noindent\item{\Bimp\ } {\smc A. Bressan}, 
On differential systems with impulsive controls, 
{\it Rend. Sem. Mat. Univ. Padova}, 
Vol. {\bf 78}, (1987), pp.  227-235.
\medskip
%
%
%
\noindent\item{\CLRSfslf\ } {\smc F.H. Clarke, Yu.S. Ledyaev, L. Rifford,
R.J. Stern,}
Feedback stabilization and Lyapunov functions,
{\it SIAM J. Control Optim.},
 {\bf 39}, (2000), no. 1, pp. 25-48.
\medskip
\noindent\item{\CLSS\ } {\smc F.H. Clarke, Yu.S. Ledyaev, E.D. Sontag,
A.I. Subbotin,}
Asymptotic controllability implies feedback stabilization,
{\it IEEE Trans. Autom. Control}, {\bf 42} (1997), pp. 1394-1407.
\medskip
%
%
\noindent\item{\CLSWnac\ } {\smc F.H. Clarke, Yu.S. Ledyaev, R.J. Stern,
P.R. Wolenski,} 
{\it Nonsmooth Analysis and Control Theory,} {\bf 178}, Springer-Verlag,    
New York, (1998).
\medskip
\noindent\item{\Riscl\ } {\smc L. Rifford}, 
Existence of Lipschitz and semiconcave control-Lyapunov
functions, 
{\it SIAM J. Control Optim.}, {\bf 39}
(2000), no. 4, pp. 1043-1064.
\medskip
\noindent\item{\Risclsf\ } {\smc L. Rifford}, 
Semiconcave control-Lyapunov
functions and stabilizing feedbacks, (2000),
preprint.
\medskip
%
%
%
%
\noindent\item{\Ss\ } {\smc E.D. Sontag}, 
Stability and stabilization: discontinuities and the
effect of disturbances, in \,
{\it Proc. NATO Advanced Study Institute - 
Nonlinear Analysis, Differential Equations, and Control}, 
(Montreal, Jul/Aug 1998), F.H. Clarke and R.J. Stern eds.,
Kluwer, (1999),\  pp. 551-598.
\medskip
%
%
%
%

\end